\documentclass[10 pt]{article}

\usepackage{amsfonts}
\usepackage{amsmath}
\usepackage{amscd}
\usepackage{amssymb}
\usepackage{latexsym}
\usepackage{graphicx}
\usepackage{amsrefs}

\newcommand{\abvr}{(\bbk,-\aak)^T}
\newcommand{\abvhr}{(-\betah_k,\alpha_k)^T}
\newcommand{\ida}{\stackrel{(a)}{=}}
\newcommand{\idb}{\stackrel{(b)}{=}}
\newcommand{\D}{{\cal D}}
\newcommand{\Do}{{{\cal D}}_1}
\newcommand{\Dt}{{{\cal D}}_2}
\newcommand{\diff}{\cl_+-\cl_-}
\newcommand{\?}{\stackrel{?}{=}}
\newcommand{\diag}{\mbox{Diag}}
\newcommand{\abv}{\bpm \bbk \\ -\aak\epm}
\newcommand{\abvj}{\bpm -\beta_j \\ \alpha_j\epm}
\newcommand{\abvh}{\bpm -\betah_k \\ \alphah_k\epm}
\newcommand{\rqv}{\bpm r_k \\ q_k \epm}
\newcommand{\aak}{\alpha_k}
\newcommand{\bbk}{\beta_k}
\newcommand{\alphah}{\hat{\alpha}}
\newcommand{\betah}{\hat{\beta}}
\newcommand{\deltah}{\hat{\delta}}
\newcommand{\gammah}{\hat{\gamma}}
\newcommand{\ggk}{\gamma_k}
\newcommand{\ddk}{\delta_k}
\newcommand{\ce}{{\cal E}}
\newcommand{\nut}{\dot{\nu}}
\newcommand{\rdt}{\dot{r}}
\newcommand{\qdt}{\dot{q}}

\newcommand{\fioh}{\hat{\phi}^{(1)}}
\newcommand{\fith}{\hat{\phi}^{(2)}}
\newcommand{\sioh}{\hat{\psi}^{(1)}}
\newcommand{\sith}{\hat{\psi}^{(2)}}
\newcommand{\sio}{\psi^{(1)}}
\newcommand{\sit}{\psi^{(2)}}
\newcommand{\fio}{\phi^{(1)}}
\newcommand{\fit}{\phi^{(2)}}

\newcommand{\phih}{\hat{\phi}}
\newcommand{\psih}{\hat{\psi}}
\newcommand{\ah}{\hat{a}}
\newcommand{\bh}{\hat{b}}

\newcommand{\ed}{E^-}
\newcommand{\eu}{E^+}

\newcommand{\srl}{\sum_{k=-\infty}^\infty}
\newcommand{\jkm}{J_k^-}
\newcommand{\jkp}{J_k^+}
\newcommand{\opi}{{\cal{I}}}
\newcommand{\om}{\omega}
\newcommand{\cc}{{\cal{C}}}
\newcommand{\ck}{{\cal{K}}}
\newcommand{\cl}{{\cal{L}}}
\newcommand{\cj}{{\cal{J}}}
\newcommand{\cR}{{\cal{R}}}
\newcommand{\cjh}{\hat{\cal{J}}}
\newcommand{\bpm}{\begin{pmatrix}}
\newcommand{\epm}{\end{pmatrix}}
\newcommand{\ov}{\overline}
\newcommand{\pf} {\noindent{\it Proof.  }}
\newcommand{\epf}{\hfill{$\square$}\medskip}
\newtheorem{cor}{Corollary}
\newtheorem{thm}{Theorem}
\newtheorem{lem}{Lemma}
\newtheorem{rmk}{Remark}
\newtheorem{prop}{Proposition}

\begin{document}

\title{A Bi-Hamiltonian Structure for the Integrable, Discrete Non-Linear Schr\"odinger System}

\author{Nicholas M. Ercolani \\ 
{\small{\it Department of Mathematics, University of Arizona, }} \\ 
{\small{\it Tucson, AZ 85721-0089, USA -- ercolani@math.arizona.edu}}\\ \\
Guadalupe I. Lozano \\  
{\small{\it Department of Mathematics, University of Michigan,  }} \\
{\small{\it Ann Arbor, MI 48109-1043, USA -- guada@umich.edu}}}

\date{}

\maketitle

\begin{abstract}

This paper shows  that the AL (Ablowitz--Ladik) hierarchy of (integrable) equations can be explicitly viewed as a hierarchy of commuting flows which: (a) are Hamiltonian with respect to both a standard, local Poisson operator $\cj$, and a new non-local, skew, almost Poisson operator $\ck$, on the appropriate space; (b) can be recursively generated from a {\it recursion} operator $\cR=\ck\cj^{-1}$.  In addition, the proof of  these facts relies upon two new pivotal resolvent identities which suggest a general method for uncovering bi-Hamiltonian structures for other families of discrete, integrable equations.

\end{abstract}
\bigskip

\noindent{\small{\bf Keywords:} discrete integrable equations, lattice dynamics, inverse scattering, Poisson geometry, bi-Hamiltonian structures.}

\section{Introduction}

The Ablowitz-Ladik (AL) system is one of the most studied discrete
integrable systems of \emph{soliton} type. It can be thought of as
an integrable discretization of the Nonlinear Schr\"odinger (NLS)
equation. However, it has recently received a great deal of
attention as the background lattice system in a variety of
modelling applications including optical fiber arrays \cite{AW, ADP}, chaos in 
dispersive numerical schemes \cite{Sc, CEM, AHS}, and linkage dynamics \cite{DS, guada}.

AL is an infinite integrable system of soliton type by which one
generally means that an infinite family of constants of motion for
the AL flow can be constructed through the inverse scattering
transform (IST) associated to a particular (discrete) eigenvalue
problem. In addition, the IST framework provides a mechanism
through which large families of solutions, such as multi-solitons
can be explicitly found.

AL is also an integrable Hamiltonian system with respect to a
natural local Poisson structure. Indeed this
structure is the discretization of the natural Poisson structure
for NLS. For many integrable partial differential equations (PDEs)
of soliton type it has been found that these systems are
Hamiltonian with respect to two \emph{distinct} Poisson
structures. Moreover, these distinct structures are compatible in
a way that enables one to construct from them a recursion operator
which generates the complete hierarchy of commuting flows which is
the signature of a completely integrable Hamiltonian system and
which characterizes its Poisson geometry. If a system has two distinct
compatible Poisson structures it is referred to as \emph{bihamiltonian} (see Section 2.) In many
specific instances the bihamiltonian  structure provides a means
to relate the Poisson geometry to the IST via resolvent relations
associated to the linear eigenvalue problem.

Examples of bihamiltonian structures related to an IST are much
rarer for discrete integrable systems. The principle goal of this
paper is to build and understand this relation for AL.

This paper is organized as follows. In Section 2 we review the
second Poisson structure for NLS which is due to Magri
\cite{Magri}.  Our approach is to explicitly re-derive Magri's structure
from a Wronskian construction of commuting flows
for soliton PDEs due to Calogero and Degasperis \cite{CalogDegas}.
We review the derivation of the second structure for the continuous case (NLS)
because it is an important guide for identifying the second
structure in the discrete case. A complete and detailed proof that
the Magri structure is Poisson may be found in \cite{guada}.

In Section 3 the necessary inverse scattering background for the
discrete case is presented including the derivation, from a
generalized Wronskian relation, of the basic operators $\cl_+$ and
$\cl_-$ in terms of which the recursion operator $\cR$ and its
properties are developed. In the last part of this section, the resolvent identities
related to these operators, which are fundamental to the main
results of this paper, are presented.  These resolvent identities represent a novel
contribution to the literature on scattering theory for discrete systems.  Full
details of their derivation may be found in \cite{guada}.

Finally in Section 4 the Poisson-geometric interpretation of the
recursion operator $\cR$ is made.  We prove that AL has a bihamiltonian character
(as defined in this section.)  To establish this it suffices to show that the
second structure, $\ck$, for AL is almost Poisson. It is an
interesting open problem to determine whether or not $\ck$ is Poisson. This
topic and other potential investigations are discussed in the conclusions.

\section{A Generalized Wronskian Approach to the Poisson Structure of NLS}

To illustrate the rationale of our approach to the geometry of the AL equations \cite{guada}, 
we first look at the parallel continuous object (the NLS hierarchy), where the geometry is already understood.

The NLS equation \eqref{NLS} is a well-known integrable PDE.  As
such, it possess, infinite families of linearly independent
constants of motion in involution, and families of explicit
special solutions such as N-soliton solutions. In the late 70's,
Magri showed that NLS has a bi-Hamiltonian nature \cite{Magri}.
This property (which now characterizes many integrable PDEs),
means that the NLS equation can be written as a Hamiltonian system
with respect to two different, independent Poisson brackets. By
composing the Poisson operators (denoted $J$ and $K$, for
instance) induced by these brackets in the appropriate manner, one
obtains a recursion operator, $R$, capable of generating a
commuting family of Hamiltonian flows, which include NLS.
Typically, one of the aforementioned Poisson operators, say $J$,
is invertible and then $R=KJ^{-1}$.

The central idea is to recover Magri's Poisson structure for NLS
from two integro-differential operators associated to the AKNS
hierarchy, a collection of commuting integrable evolution
equations which includes NLS.

The first of these operators, which we denote $L_+$, was
constructed by Calogero and Degasperis using a generalized
Wronskian technique \cite{CalogDegas}.  The basic idea behind the
generalized Wronskian approach to integrable evolution equations,
is to generate a set of \emph{scattering} relations between the
asymptotic behaviors (as $x\to +\infty$ and $x\to -\infty$) of the
wave-function solutions of some initial eigenvalue problem, such
as the Zakharov--Shabat problem
\begin{equation}
\bpm {\psi_1}_x \\ {\psi_2}_x \epm =\bpm -ik & q \\ r & ik \epm
\bpm \psi_1 \\ \psi_2 \epm.\label{ZS-cont}
\end{equation}

The second operator, $L_-$,  appearing in Calogero--Degasperis' work and giving rise to the AKNS hierarchy is related to $L_+$ in a sort of ``adjoint" way, as described below.

The specific form of the operators $L_+$ and $L_-$ (which act on a certain space $C$  of rapidly decaying complex-valued functions $(q(x), r(x))^T$--the potentials,) suggests the construction of an
anti-symmetric operator $K$ arising from $L_+, L_-$.

As shown in \cite{guada} (see below) one can indeed use $L_+$ and $L_-$ to define two
geometrically meaningful operators $K$ and $R$ which act on $C$.
The first one is Poisson, and gives rise to Magri's bracket for
NLS.  The second one is a recursion operator for the AKNS
hierarchy which may in fact be portrayed as the composition of $K$
with the standard Poison structure for AKNS (given essentially as
multiplication by $i$).

The specific structure and relationship between $L_+$ and $L_-$ in
the continuous setting will eventually guide the construction of
their discrete counterparts $\cl_+$ and $\cl_-$ and also motivate
the definition of $\ck$ and $\cR$ --discrete analogs of $K$ and
$R$ for the AL hierarchy \cite{guada}.  These results are presented 
next.

\subsection{Obtaining Magri's Poisson structure, $K$, for NLS through $L_+$ and $L_-$}

In their 1976 paper, Calogero and Degasperis begin with the the eigenvalue problem \eqref{ZS-cont} and use a generalized Wronskian technique to arrive at the following formally defined class of
integrable equations:
\begin{equation}
(r_t(x,t), -q_t(x,t))^T=\gamma(L_+)(r(x,t),q(x,t))^T, \label{CalogDegasEvol}
\end{equation}
where $\gamma$ is an entire function of the integro-differential
operator
\begin{equation}
L_+=\frac{1}{2i} \left[ \bpm 1 & \ 0 \\ 0 & -1 \epm
\partial_x + 2 \bpm  {rI_+q}  &  {-rI_+r} \\ {qI_+q}  &  {-qI_+r}
\epm \right], \ \ I_+=\int_x^{+\infty}(\ \cdot\ ) dy .\label{contL+}
\end{equation}

Perhaps one of the best-known procedures for generating
hierarchies of non-linear integrable evolution equations is due to
the work of Ablowitz, Kaup, Newell and Segur {\bf(AKNS)}
\cite{AKNS}. In our context of interest, their method utilizes
once again the Zakharov--Shabat eigenvalue problem together with a
second linear operator prescribing the time evolution of the
wavefunctions. The hierarchy of integrable evolution equations
arises then as a series of compatibility conditions associated to
the linear problems just described.  All the evolution equations
in the hierarchy stem from the same eigenvalue problem, yet each
corresponds to a different time-evolution operator.

Calogero and Degasperis write the AKNS hierarchy as
\begin{equation}
(r_t(x,t),-q_t(x,t))^T=-2A(L_-)(r(x,t),q(x,t))^T, \label{AKNS}
\end{equation}
where
\begin{equation}
L_-=\frac{1}{2i} \left[ \bpm 1 & \ 0 \\ 0 & -1 \epm
\partial_x - 2 \bpm  {rI_-q}  & {-rI_-r} \\ {qI_-q}  &  {-qI_-r}
\epm \right], \ \ I_-=\int_{-\infty}^x (\ \cdot\ ) dy. \label{contL-}
\end{equation}
($A$ is an arbitrary entire function), and show that \eqref{CalogDegasEvol} and \eqref{AKNS} are
equivalent.

The proof of this fact is entails showing that
\begin{equation}
L_+^n((r,q)^T)=L_-^n((r,q)^T), \label{L+hier=L-hier}
\end{equation}
which follows, in turn, from the observation that $L_+^n((r,q)^T)$
is in the kernel of the difference operator $L_--L_+$ --the
continuous ``kernel condition." This can be verified directly for
$n=0$ and other low-values of $n$. The proof for general $n$
(which relies on the definition of two auxiliary functions of the
complex variable $z$) can be found in the Appendix of the original
Calogero--Degasperis paper.

We remark that, under the identification $r=-\ov q$, the system of coupled integrable evolution equations \eqref{AKNS}
reduces to a system of evolution equations for a single field $q$, comprising the NLS hierarchy, which contains the well-known NLS equation
\begin{equation}
-iq_t=q_{xx}+2|q|^2q. \label{NLS}
\end{equation}
In other words, in the reduction $r=-\ov q$, \eqref{AKNS} becomes the NLS family of equations.

Let us now consider the space of complex, vector valued functions of a real variable $x$ given by
$$C=\left\{ (q(x), r(x))^T \in {\mathbb C}^2 :
|q|,|r|\rightarrow 0 \mbox{ as } x \rightarrow\pm\infty\right\}.
$$
The tangent bundle to $C$ may be endowed with a non-degenerate bilinear form locally described by the inner product
\begin{equation}
\langle ( v_1, v_2)^T, ( w_1, w_2)^T\rangle
 = - \int v_1 w_2 + v_2 w_1 \ \ dx,\label{contInnP}
\end{equation}
where, $(v_1, v_2)^T, (w_1, w_2)^T$ are tangent vectors to $C$ at
the point $(q,r)^T$, and $\int f\ dx=\int_{-\infty}^{\infty} f \
dx$. One may also define a skew-adjoint operator $J=\diag(-i, i),$
on the tangent bundle to $C$ which, together with the inner
product just defined, gives rise to the complexified standard
Poisson bracket on $C$, that is,
$$
\{F,G \}_J = \left\langle \nabla F, \diag(-i,i) \nabla G \right\rangle = \left\langle ( -\delta_r F,
-\delta_q F)^T, \diag(-i,i) (-\delta_r G,
-\delta_q G)^T\right\rangle.
$$
Here, $F$ and $G$ are functionals on $C$ which become real-valued in the reduction $r_k=\overline{q_k}$ defining the standard setting for NLS.  This reality condition implies that $F$ and $G$ must be symmetric (or anti-symmetric) in $r$ and $q$. Furthermore, $(\delta_q (\ \cdot\ ), \delta_r (\ \cdot\ ))^T$ denotes the variational derivative, whereas $\nabla (\ \cdot\ )=(-\delta_r (\ \cdot\ ), -\delta_q (\ \cdot\ ))^T$ gives the functional gradient of each such functional, with respect to the given inner product.

We now re-consider the operators $L_+$ and $L_-$  in the context of the function space $C$.

The presence of the integral operators $I_+$ and $I_-$ suggests that $L_+$ and $L_-$ may be combined so as to yield an \emph{anti-symmetric} operator, $K$ related to the hierarchy. Indeed, as we will show below, the special properties of $L_-+L_+$ and $L_--L_+$ allows us to recover both: the (second) Poisson structure for NLS discovered by Magri \cite{Magri} (prescribed by the Poisson operator $K$), and the recursion operator $R=KJ^{-1}$ associated to the hierarchy of flows, in the current setting.
We begin by defining,
\begin{equation}
K \doteq \bpm 1 & 0 \\ 0 & 1\epm\partial_x + \bpm  {-q(I_--I_+)r} &  {q(I_--I_+)q} \\ {r(I_--I_+)r}  &  {-r(I_--I_+)q}\epm. \label{Kcont}
\end{equation}
This operator $K$, just as $J$, acts on the tangent space to $C$ at a point and defines a (point-dependent) operator on the tangent bundle of $C$.
\begin{prop} \label{PROP:KcontPoisson}
The bracket $\{F,G \}_K = \left\langle\nabla F, K \nabla G
\right\rangle$ associated to $K$ is skew-symmetric and satisfies
the Jacobi identity.  Hence, $K$ is Poisson. \hfill $\square$
\end{prop}
A proof by direct-calculation of Proposition
\ref{PROP:KcontPoisson}, can be found in \cite{guada}.  The
technique relies on a standard integration by parts formula (see
Appendix A in \cite{guada}, for instance.)

There are at least two possible avenues to establish a precise
link between the Poisson operator $K$ and the hierarchy of
integrable evolution equations \eqref{AKNS}. The first one begins
with the observation that, in the reduction $r=-\ov q$, Hamilton's
equations associated to $K$ become those given by Magri's Poisson
bracket \cite{Magri}. Since the integrable hierarchy in the
previous section is the NLS hierarchy for $r=-\ov q$, and Magri's
bracket is known to prescribe a bi-Hamiltonian structure for NLS
(together with the Poisson bracket given as multiplication by
$-i$), the above observation connects $K$ with the AKNS hierarchy
via the identification $r=-\ov q$.

A second and perhaps more constructive approach relies on directly
checking that $K$ and $J$ form a compatible pair of Poisson
structures and that the recursion operator $KJ^{-1}$ does indeed
generate the AKNS hierarchy of flows.  This approach has the
advantage of giving an explicit description of the bi-hamiltonian
nature of the evolution equations in the AKNS hierarchy.

Recall that the compatibility of two Poisson operators such as $K$ and $J$ amounts to proving that their sum is also a Poisson operator on the manifold in question.  Two compatible Poisson operators give rise to a bi-Hamiltonian system when there exists a vector field that is Hamiltonian with respect to both Poisson structures \cite{Magri}.  In this case, one can generate a hierarchy of bi-Hamiltonian vector fields, by recursively applying the composition of one of the Poisson operators with the inverse of the other to the original bi-Hamiltonian field (assuming, of course, that one of the Poisson operators is indeed invertible.)  We define
\begin{align}
R = KJ^{-1}
  & =  \left[ \bpm 1 & 0 \\ 0 & 1\epm \partial_x
+ \bpm {-q(I_--I_+)r}
       &  {q(I_--I_+)q} \\ {r(I_--I_+)r}  & {-r(I_--I_+)q}\epm \right]
 \bpm i & 0 \\ 0 & -i\epm \nonumber \\
       & =  \bpm 0 & 1 \\ -1 & 0\epm (L_-+L_+) \bpm 0 & -1 \\ 1 &
       0\epm, \label{RecCont}
\end{align}
and consider the functional $H=-\int {rq}$.  Then, $J\nabla H = \diag(-i,i) (q,r)^T=(-iq,ir)^T,$ and $X_n = R^n J\nabla H$ defines a hierarchy of bi-Hamiltonian vector fields. Hamilton's equations associated to these fields are
\begin{equation}
((q,r)^T)_t = R^n J\nabla H. \label{Ham}
\end{equation}
Notice that $X_0 = J\nabla H=(-iq,ir)^T$ and $X_1 = K\nabla H=(q_x,r_x)^T$.  The standard NLS flow is
prescribed by
\begin{equation}
X_3=i\left[( q_{xx}, -r_{xx})^T + (-2rq^2,
2qr^2)^T \right]. \label{NLSrq}
\end{equation}
The hierarchy \eqref{Ham} is the NLS Hierarchy.

We shall next show that the NLS hierarchy \eqref{Ham} is equivalent to a sub-hierarchy of the family of Calogero--Degasperis flows \eqref{CalogDegasEvol}.  In fact, we will demonstrate that the entire Calogero--Degasperis hierarchy can be generated from polynomials in $R$.  The argument's basic idea is to write the recursion operator $R$ in terms of $L_-+L_+$ and exploit the relationship between $L_+$ and $L_-$ given in \eqref{L+hier=L-hier}.
\begin{prop}\label{PROP:contL+hier=contRhier}
The NLS Hierarchy, obtained by recursively applying the recursion operator $R$ to the (Hamiltonian) field $(-iq,ir)^T$, is equivalent to a sub-hierarchy of the Calogero--Degasperis family
\eqref{CalogDegasEvol}.
\end{prop}

\pf  We re-write Hamilton's equations \eqref{Ham}, using the form of the recursion operator given in \eqref{RecCont}, to obtain
\begin{equation}
((q,r)^T)_t= \sigma_2 (L_-+L_+)^n(r,q)^T, \label{Ham0}
\end{equation}
where
\begin{equation}
\sigma_2=\bpm 0 & -i\\i & 0\epm
\end{equation}
is the standard Pauli matrix.

Using induction on $n$ together with formula \eqref{L+hier=L-hier} one sees that
\begin{equation}
(L_-+L_+)^n(r, q)^T=2^n(L_+)^n(r,q)^T, \label{Ham1}
\end{equation}
so that \eqref{Ham0} becomes
\begin{equation}
((q,r)^T)_t = 2^n \sigma_2 L_+^n (r,q)^T. \label{Ham2}
\end{equation}
Multiplying both sides of \eqref{Ham2} by $-i\sigma_2$, we obtain
\begin{equation}
 ((-r,q)^T)_t =-2^n i L_+^n (r,q)^T, \label{Ham3}
\end{equation}
which is clearly a sub-hierarchy of the Calogero--Degasperis
family of flows. \epf

Using Weierstrass approximation one can replace the monomial in \eqref{Ham3} by an arbitrary analytic function, $A$, of $L_+$ and finally re-write the equation as $((r, -q)^T)_t +A(L_+) (r,q)^T= 0,$ which is the form of the flow derived by Calogero and Degasperis.  Hence,
\begin{cor}\label{contL+fam=contRfam}
The entire Calogero--Degasperis family of flows may be formally generated from polynomials in $R$.\hfill $\square$
\end{cor}

\section{Inverse Scattering Preliminaries}

This section begins with the direct scattering formulation for the Zakharov--Shabat eigenvalue problem \eqref{ZS-disc} originally studied by Ablowitz and Ladik in connection with the integrable
equations bearing their name.

After identifying (and characterizing) two pairs of eigenfunctions
$\phi_k(z)$, $\phih_k(z)$, $\psi_k(z)$, $\psih_k(z)$
(which behave nicely at either $-\infty$ or $\infty$),
we focus on a Wronskian relation introduced by Ladik and Chiu as a means to
study a family of integrable evolution equations arising from the
linear problem \eqref{ZS-disc} (which includes AL.)

The generalized Wronskian identity presented by Ladik and Chiu
provides, once again, a way to relate the evolution of the
potentials $r_k$ and $q_k$ (defining \eqref{ZS-disc}), and the
kinematics of the scattering coefficients (i.e., the parameters
specifying the asymptotic behavior of the eigenfunctions mentioned
above.) By so doing, it also singles out two sum-difference
operators $L$ and $L^{-1}$, which together generate the collection
of (integrable) equations prescribing the potentials' evolution.

As we will soon see, the specific form and relationship between
$L$ and $L^{-1}$ will guide the construction of two new operators
$\cl_+$ and $\cl_-$, associated to a sub-hierarchy of the
Chiu--Ladik equations which encompasses the AL system \cite{guada}.  The latter
operators can be viewed as discrete analogs of $L_+$ and $L_-$ in
more than one way.  In particular, they will play a key role in
unveiling the Poisson-geometric picture behind AL, as we begin to
indicate below.

Exploiting the time-independence of one of the scattering coefficients (an essential feature of these kind of integrable equations), one can obtain both: a sequence of constants of motion, and a generating function for the associated hierarchy of variation (gradient) fields.  It turns out that it is also possible to explicitly write down  (resolvent-type) identities giving the infinite hierarchy of gradient fields in terms of
powers of $L$ and $L^{-1}$.  As explained below, these identities will be instrumental for exhibiting the bi-Hamiltonian nature of the AL hierarchy in Section \ref{geom}.

\subsection{A generalized Wronskian leading to AL flows}

In 1976, Chiu and Ladik orchestrated a generalized Wronskian identity based on the discrete version of
the Zakharov--Shabat eigenvalue problem
\begin{equation}
\bpm \nu_1 \\ \nu_2 \epm_{k+1} =\bpm z & q_k \\ r_k & 1/z \epm
\bpm \nu_1 \\ \nu_2 \epm_{k}=\ce_k \bpm \nu_1 \\ \nu_2 \epm_{k},
\label{ZS-disc}
\end{equation}
which generated a large family of integrable evolution equations
associated to AL. Here, $z\in{\mathbb C}$ is the eigenvalue
parameter, and the complex-valued potentials $q_k$ and $r_k$ are
assumed to vanish rapidly for $k\rightarrow\pm\infty$.

Let $\Phi_k(z)$ and $\Phi'_k(z)$ respectively denote matrix solutions to \eqref{ZS-disc} and to a \emph{second} eigenvalue problem of the same form but for a new set of potentials $q'_k$, $r'_k$. The initial Wronskian identity proposed by Chiu and Ladik relates $\Phi_k$ and $\Phi'_k$ to $a, \ah, b, \bh$ and $a', \ah', b', \bh'$ (the ``unprimed" and ``primed" scattering parameters) via the right-hand-side of the equation below \cite{LC}:
\begin{align}
&\sum_{k=-\infty}^{\infty} P(k+1)\ {\Phi'_k}^T\  M(k,z)\  \Phi_{k+1}={\Phi'_k}^T P (k)F(k)\Phi_k|_{-\infty}^{\infty}, \nonumber \\
&\mbox{ where }\ \ M(k,z)=\left[\bpm z & r_k'\\ q_k' & 1/z \epm F(k+1)-F(k)\bpm 1/z & -q_k \\ -r_k & z \epm\right].
\label{wronsOrig}
\end{align}
Above, $P(k)=\prod_{j=k}^\infty(1-r_kq_k)$ and $F(k,z)$ is an arbitrarily chosen matrix.  The validity of \eqref{wronsOrig} results from the telescoping of the terms in the specified series, which can
be checked directly (see \cite{guada} for details.)

Due to the asymptotic decay of the potentials, the right-hand-side of \eqref{wronsOrig} makes sense and may be written in terms of the scattering parameters (for the ``primed" and ``unprimed" linear problems.)  Its precise form depends both on the choice of matrix $F(k,z)$ and on the particular eigenfunctions making up the columns of $\Phi_k(z)$.
\begin{rmk} We will use the notation $f_k$ and $f(k)$ interchangeably, depending on emphasis, to denote functions of the discrete variable $k$.  The $z$ dependence of objects such as $\Phi_k$, $F(k)$ and the scattering parameters $a, \ah$, $b, \bh$, may or may not be explicitly noted, depending on the context.
\end{rmk}
Following the work of Chiu--Ladik, we now specify the precise form of the right-hand-side of \eqref{wronsOrig} by choosing a particular fundamental matrix $\Phi_k(z)$ and determining its asymptotic form.   We begin by considering the four special solutions
\begin{align}
\phi_k(z)\sim z^k(1,0)^T, \ \ \ \phih_k(z)\sim
z^{-k}(0,1)^T, & \mbox{ for } k\rightarrow-\infty \label{thePhis} \\
\psi_k(z)\sim z^{-k}(0,1)^T, \ \ \ \ \psih_k(z)\sim
z^k(1,0)^T, & \mbox{ for } k\rightarrow\infty. \label{thePsis}
\end{align}
of \eqref{ZS-disc} and the scattering relations prescribed by
\begin{equation}
\bpm \phi_k(z) \\ \phih_k(z)\epm = \bpm b(z) & a(z)\\ \ah(z) &
\bh(z)\epm \bpm \psi_k(z) \\ \psih_k(z)\epm = S(z) \bpm \psi_k(z) \\
\psih_k(z)\epm. \label{phiScattRel}
\end{equation}
We then set $\Phi_k(z)=(\phi_k(z),\psi_k(z))$ and determine that:
\begin{equation}
\Phi_k\sim\bpm z^k & \frac{-\bh}{C_0}z^k \\ 0 & \frac{a}{C_0}z^{-k}\epm, \mbox{ for } k\rightarrow-\infty, \ \ \  \ \Phi_k\sim\bpm az^k & 0 \\ bz^{-k} & z^{-k}\epm, \mbox{ for } k\rightarrow\infty, \label{PhiAsym}
\end{equation}
where $\det S(z)=a\ah-b\bh=C_0$,  and $C_0=\lim_{k\rightarrow-\infty}P(k)=\Pi_{-\infty}^{\infty}(1-r_kq_k)$, as shown by Ablowitz, Prinari and Trubatch \cite{APT}.  (See also \cite{guada} for more details.)

We now adapt the results established by Chiu and Ladik  to our choice of $\Phi(z)$ and arrive at two Wronskian-type identities based on \eqref{wronsOrig}.   They are:
\begin{multline}
\sum_{k=-\infty}^\infty P(k+1){\Phi'_k}^T\left[ \opi \Lambda^l\bpm H_1^{(0)}(k)\\ H_4^{(0)}(k)\epm\right]\Phi_{k+1}\hspace{+2.5in} \\
\hspace{+.7in} = z^{2l}\bpm a'b & a'-a \\ 0 &
\displaystyle{\frac{a\bh'}{C_0'}}\epm + \sum_{j=1}^l
z^{2(l-j)}({\Phi'_k}^T P(k)F^{(j)}(k)\Phi_k)|_{-\infty}^\infty,
\label{wronZ2Ml}
\end{multline}\vspace{-.15in}
\begin{multline}
\sum_{k=-\infty}^\infty P(k+1){\Phi'_k}^T\left[\opi{(\Lambda^{-1})}^l\bpm h_1^{(0)}(k)\\
h_4^{(0)}(k)\epm \right]\Phi_{k+1} \\
=\frac{1}{z^{2l}}\bpm a b' & 0 \\ a -
\displaystyle{\frac{a'C_0}{C_0'}} &
\displaystyle{\frac{a'\bh}{C_0'}} \epm + \sum_{j=1}^l
\frac{1}{z^{2(l-j)}}({\Phi'_k}^T
P(k)F^{(j)}(k)\Phi_k)|_{-\infty}^\infty , \label{wronZinv2Ml}
\end{multline}
where the integro-differential operators $\Lambda$, $\Lambda^{-1}$ are given by
\begin{align}
\Lambda = \bpm E^- & 0 \\ 0 & E^+ \epm & +
\bpm {-r_k\cj_k(q_j E^-)} & {r_k\cj_k(r'_j E^+)} \\
{-q'_k\cj_k(q_j E^-)} & {q'_k\cj_k(r'_j E^+)}\epm + \bpm 0 & 0 \\
q'_k q_k E^- & -q'_k r'_k E^+\epm \nonumber \\
& +  (1-r_k q_k)\bpm {-E^-(r'_k)\cjh_k(q'_j)} &
{E^-(r'_k)\cjh_k(r_j)} \\ 0 & 0 \epm \nonumber \\
& + (1-r'_k q'_k) \bpm 0 & 0 \\ {-E^+(q_k\cjh_k(q'_j))} &
{E^+(q_k\cjh_k(r_j))}\epm , \label{lambda}
\end{align}\vspace{-.1in}
\begin{align}
\Lambda^{-1} = \bpm E^+ & 0 \\ 0 & E^-\epm & +
\bpm {r'_k\cj_k(q'_j E^+)} & {-r'_k\cj_k(r_j E^-)} \\
{q_k\cj_k(q'_j E^+)} & {-q_k\cj_k(r_j E^-)}\epm + \bpm -r'_k q'_k
E^+ & r'_k r_k E^- \\ 0 & 0\epm \nonumber \\
& +   (1-r_k q_k)\bpm 0 & 0 \\ {E^-(q'_k)\cjh_k(q_j)} &
{-E^-(q'_k)\cjh_k(r'_j)}\epm \nonumber \\
& +  (1-r'_k q'_k)\bpm {E^+(r_k\cjh_k(q_j))} &
{-E^+(r_k\cjh_k(r'_j))} \\ 0 & 0\epm, \label{lambdatilde}
\end{align}
as demonstrated in the original Chiu--Ladik work \cite{LC}.   The shift operators $E^\pm$ act in the standard fashion, $E^\pm v_k=v_{k\pm1}$, whereas $\cj_k (v(j))\ u \doteq  \sum_{j=k}^\infty v(j)u(j)$, and $\cjh_k (v(j))\ u \doteq p(k)/2+\sum_{j=k}^\infty\ [(\prod_{i=k}^{j-1}(1-r'_i q'_i)/(1-r_i q_i))(v(j)u(j))/(1-r_j q_j)]$, for $p(k)\doteq p \prod_{j=k} ^\infty(1-r'_k q'_k)/(1-r_k q_k)$ and $p$ an arbitrary constant.

As explained both in \cite{guada} and in Chiu--Ladik's paper, Eq. \eqref{wronZ2Ml} arises from an iteration in powers of $z^2$ prescribed essentially through a suitably chosen ansatz for the $M(k,z)$ factor in the summand of \eqref{wronsOrig}.  Specifically,
\begin{equation}
M^{(l)}(k)=\bpm {-z^2H_1^{(l-1)}(k)+H_1^{(l)}(k)}& 0 \\ 0 &
{-z^2H_4^{(l-1)}(k)+H_4^{(l)}(k)}\epm\ , l\geq 0 \label{Ml}
\end{equation}\vspace{-.15in}
where
\begin{equation}
M^{(0)}(k)=\bpm r_k & 0 \\ 0 & q'_k\epm=\bpm H_1^{(0)}(k) & 0 \\ 0 & H_4^{(0)}(k)\epm, \mbox{ and } \bpm H_1^{(l+1)}(k) \\ H_4^{(l+1)}(k)\epm = \Lambda \bpm
H_1^{(l)}(k) \\ H_4^{(l)}(k)\epm.
\end{equation}
We remark that such an ansatz completely determines the structure of $F(k,z)$.  The second equation emerges in a similar manner from an iteration in powers of $1/z^2$.

The pivotal identity for obtaining the precise time-evolution of
the vector potential $(r_k, q_k)$ and the scattering parameters
associated to \eqref{ZS-disc}, arises by taking linear
combinations of \eqref{wronZ2Ml} and \eqref{wronZinv2Ml} for
different values of $l>0$ \cite{guada}.  When temporal
considerations are brought to the forefront, the aforementioned
combination of \eqref{wronZ2Ml} and \eqref{wronZinv2Ml} becomes
time-dependent (through the $t$-dependent vector potential,
scattering data, etc.)  One then exploits this time dependence by
defining the scattering parameters in the ``primed" eigenvalue
problem as time-evolved variants of those in the ``non-primed"
scenario, namely,
\begin{align}
 a' &\doteq  a(t), \ \ \ \ah' \doteq \ah(t), \ b' \doteq b(t), \ \bh' \doteq \bh(t), \nonumber \\
 r'_k  &\doteq  r_k(t), \ \ q'_k \doteq q_k(t), \ \Phi'_k  \doteq  \Phi_k(t),\ C_0' \doteq C_0(t).
\label{tVariants}
\end{align}
As detailed in Chapter 3 of \cite{guada}, evaluating the right difference-quotients and limits --based on the substitutions \eqref{tVariants}-- leads to the operators
\begin{align}
L = \bpm E^- & 0 \\ 0 & E^+\epm & +  (1-r_k q_k)\bpm -p E^-(r_k) \\
-\check{p} E^+(q_k) \epm \nonumber \\ & +  (1-r_kq_k)\bpm
{-E^-\left(r_k \jkp\left(\frac{q_j}{1-r_jq_j}\right)\right)} &
 {E^-\left(r_k \jkp\left(\frac{r_j}{1-r_jq_j}\right)\right)} \\
{-E^+\left(q_k \jkp\left(\frac{q_j}{1-r_jq_j}\right)\right)} &
 {E^+\left(q_k \jkp\left(\frac{r_j}{1-r_jq_j}\right)\right)}\epm
 \nonumber \\ & +  (1-r_kq_k)\bpm
 {E^-\left(\frac{r_kq_k}{1-r_kq_k}\right)} &
{-E^-\left(\frac{r_k^2}{1-r_kq_k}\right)} \\ 0 & 0\epm\nonumber
\\ & +  \bpm -r_k \jkp(q_jE^-) & r_k \jkp(r_jE^+) \\
-q_k \jkp(q_jE^-) & q_k \jkp(r_jE^+)\epm + \bpm 0 & 0 \\ q_k^2 E^- & -q_kr_k E^+\epm, \label{L}
\end{align}
\begin{align}
L^{-1} = \bpm E^+ & 0 \\ 0 & E^-\epm & +  (1-r_k q_k)\bpm p E^+(r_k) \\
\check{p} E^-(q_k)\epm \nonumber \\ & + (1-r_kq_k)\bpm
 {E^+\left(r_k \jkp\left(\frac{q_j}{1-r_jq_j}\right)\right)} &
{-E^+\left(r_k \jkp\left(\frac{r_j}{1-r_jq_j}\right)\right)} \\
 {E^-\left(q_k \jkp\left(\frac{q_j}{1-r_jq_j}\right)\right)} &
{-E^-\left(q_k \jkp\left(\frac{r_j}{1-r_jq_j}\right)\right)}\epm
 \nonumber \\ & +  (1-r_kq_k)\bpm 0 & 0 \\
{-E^-\left(\frac{q_k^2}{1-r_kq_k}\right)} &
 {E^-\left(\frac{q_kr_k}{1-r_kq_k}\right)}\epm\nonumber
\\ & +  \bpm r_k \jkp(q_jE^+) & -r_k \jkp(r_jE^-) \\
q_k \jkp(q_jE^+) & -q_k \jkp(r_jE^-)\epm + \bpm -r_kq_k E^+ & r_k^2 E^- \\ 0 & 0\epm, \label{Linv}
\end{align}
(arising from $\Lambda$ and $\Lambda^{-1}$ respectively), and to the following (simple!) time-evolution equations for the scattering parameters:
\begin{equation}
((a(z),\ah(z))^T)_t= 0, \ \ \ \  ((-b(z),\bh(z))^T)_t = \om(z^2)(b(z),\bh(z))^T, \ \ \ \ (C_0)_t = 0. \label{ab-evol}
\end{equation}
\begin{rmk} The operator $\jkp$ present in \eqref{L} and \eqref{Linv} is defined by $\jkp (u_j)\doteq \sum_{j=k}^\infty u_j$ and is the discrete analog of the integral operator $I_+$, appearing in the expression for the continuous integro-differential operator $L_+$ given in \eqref{contL+}.  The symbols
$p,\check{p}$ denote (discrete) integration constants.
\end{rmk}
Obtaining evolution equations for the potentials $r_k(t)$,
$q_k(t)$, requires more careful arguments which exploit the
analytic properties of the Jost functions associated to the linear
problem \eqref{ZS-disc} and its four special solutions
\eqref{thePhis}--\eqref{thePsis}.  Through the proof of Theorem 1
in \cite{guada}, one arrives at the following time-flows:
\begin{equation}
((-r_k, q_k)^T)_t = \om(L)(r_k,q_k)^T, \label{PotEvol}
\end{equation}
where $\om(x)\doteq\om_1(x)+\om_2(x^{-1})$ and $L$ and $L^{-1}$ are the sum-difference operators above.

The family of discrete integrable evolution
equations\eqref{PotEvol} arising from the Chiu--Ladik Wronskian
approach includes several discrete versions of mKdV equations as
well as the well-known discretization of NLS constructed by
Ablowitz and Ladik, namely,
\begin{equation}
i(q_k)_t=(1+|q_k|^2)(q_{k+1}+q_{k-1})\label{ALqqbar},
\end{equation}
(up to the linear term $-2q_k$.)  The latter --which we will also call the AL equation-- can essentially be obtained from one of the simplest symmetric $\om(L)$, namely $L+L^{-1}$. Indeed, as one can directly check,
\begin{equation}
((-r_k,q_k)^T)_t = -i(L+L^{-1}) (r_k,q_k)^T = -i(1-r_kq_k)(E^++E^-)(r_k,q_k)^T
\label{ALrq}
\end{equation}
In the reduction $r_k=-\overline{q}_k$, Eq. \eqref{ALrq} becomes
\eqref{ALqqbar}.

\subsection{The constant of motion $\log a$ and associated gradient fields}

Recall that the scattering parameter $a$ is an analytic function
of the complex variable $z$ on the complement of the unit disc,
and can therefore be represented by its Laurent series expansion
(see \cite{guada} and \cite{APT}, for more details.) Since $a$ is
time-independent --as made evident through the evolution equations
\eqref{ab-evol}-- so are the coefficients of its Laurent series.
Such coefficients then constitute an infinite family of constants
of motion associated with the flows \eqref{PotEvol}. The same
reasoning shows that the coefficients of the Taylor series for
$\ah$ (on $|z|<1$) also give an infinite family of integrals.

In their 1975 paper, Ablowitz and Ladik \cite{AL2} utilize a recursive technique to arrive at the
collection of constants of motion arising from $\log\ah$:
\begin{equation}
\hat{C}_1=\sum q_{k}r_{k+1}, \ \ \  \hat{C}_2=\sum r_{k+1}q_{k-1}(1-r_k q_k)-\frac{1}{2}q_k^2r_{k+1}^2, \ \ \
\cdots \label{ooooo}
\end{equation}
The same formal method (outlined in \cite{guada}), may be used to compute the hierarchy of
integrals associated to $\log a$, namely,
\begin{equation}
C_1=\sum r_{k}q_{k+1}, \ \ \ C_2=\sum\ q_{k+1}r_{k-1}(1-r_k q_k)-\frac{1}{2}r_k^2q_{k+1}^2, \ \ \
\cdots \label{oooo}
\end{equation}
One may also obtain a hierarchy of ``symmetric" constants of motion, by combining the above as
$C_0, \ \ C_1+\hat{C}_1, \ \ C_2+\hat{C}_2, \ \ \cdots,$ where $C_0=\prod_{k=-\infty}^\infty(1-r_kq_k)$. (See \eqref{PhiAsym}, \eqref{ab-evol}.)

We now determine  a hierarchy of gradient fields associated to the conserved quantities arising from the series expansion of $\log a$ and $\log\ah$ just discussed.  Our approach here draws on the standard treatment of the parallel derivation for the modified Korteweg-DeVries equation (mKdV), a well-known integrable PDE on the real line \cite{newell}.

The argument is based upon considering a variation $(\qdt_k,\rdt_k)^T=(a_k, b_k)^T$ of  $(q,r)^T$ along the space of rapidly vanishing (vector) potentials, while keeping the eigenvalue $z$ fixed.  By finding the expression for the induced variation on the special eigenfunctions $\phi_k, \phih_k, \psi_k, \psih_k$ (the special solutions given in \eqref{thePhis}, \eqref{thePsis} to the linear problem \eqref{ZS-disc}), one is able to precisely calculate the effect of the potentials' perturbation on the scattering parameters $a, \ah,
b, \bh$.  For instance,
\begin{equation}\vspace{-.2in}
\dot{\psi}_n=-\Phi(n)\sum_{k=n}^\infty \Phi(k+1)^{-1}\bpm 0 & \qdt_k\\ \rdt_k & 0 \epm \ \psi_k,\label{psiVariation}
\end{equation}
leads to
\begin{equation}
\frac{C_0}{a}\dot{\left(\frac{a}{C_0}\right)}=\dot{\left[\log a-\sum\log(1-r_kq_k)\right]}=\sum \frac{P(k+1)}{a}\left(-\fio_{k+1}\sio_k\ \rdt_k+\fit_{k+1}\sit_k\
\qdt_k\right), \label{aCoFinalVariation}
\end{equation}
in the manner outlined in Appendix \ref{APP:I} of this work. This,
in turn, allows us to deduce an explicit expression for the
variational derivative $\delta_k \log a\doteq (\delta \log
a/\delta q_k,\delta \log a/\delta r_k)$ of $\log a$, namely,
\begin{equation}
\delta_k\left(\log a-\sum\log (1-r_kq_k)\right)=\frac{P(k+1)}{a}\big(\fit_{k+1}\sit_k,-\fio_{k+1}\sio_k\big)^T.\label{yeah!}
\end{equation}
By Laurent-expanding the expression for this generating function one can then obtain the desired hierarchy of ``gradient" fields, each given by the variational derivative of the corresponding coefficient in the expansion of $\log a$.

The results just outlined, together with the supporting arguments presented in Appendix \ref{APP:I}, prove the
following theorem.
\begin{thm} \label{THM:grads}  The (z-dependent) gradients of $\log a(z)$ and $\log\ah(z)$ are generating functions for the hierarchy of variation fields $\delta_kC_n$ and $\delta_k\hat{C}_n$ associated to the constants of motion \eqref{oooo} and \eqref{ooooo} respectively. Their explicit expressions are:
\begin{align}
\delta_k\log a+\delta_k H_0&=\frac{P(k+1)}{a} (\fit_{k+1}\sit_k,\  -\fio_{k+1}\sio_k)^T, \label{aGrads-psi2} \\
\delta_k\log \ah + \delta_k H_0 &=\frac{P(k+1)}{\ah}(-\fith_{k+1}\sith_k,\ \fioh_{k+1}\sioh_k)^T,\label{ahGrads-psi1}
\end{align}
where $H_0=-\sum_{k=-\infty}^\infty\log (1-r_kq_k)=-\log C_0$ and $P(k+1)=\prod_{j=k+1}^\infty (1-r_kq_k)$.\hfill$\square$
\end{thm}
Notice that, on the unit circle $|z|=1$, both Eq. \eqref{aGrads-psi2} and Eq. \eqref{ahGrads-psi1} make sense (as the scattering coefficients $a(z),\ah(z)$ are holomorphic on the regions $|z|>1$ and $|z|<1$ of ${\mathbb C}$ respectively.) One could then combine them to obtain a ``symmetric" representation of the variational derivatives.  This line of reasoning leads to the symmetric constants of motion stemming from the combination of \eqref{oooo} and \eqref{ooooo}.

\subsection{Constructing the generating operators: $\cl_+$, $\cl_-$ and $\cR$}

Prima facie, one observes several qualitative analogies between the operators $L_+^2$ (defined in Eq. \eqref{contL+}) associated to continuous NLS and $L+L^{-1}$.  Both are second order operators, and both may be obtained via parallel (generalized Wronskian) techniques (in the discrete
and continuous settings respectively.) Furthermore, they each give rise to (continuous and discrete) integrable versions of NLS in a parallel manner:
\begin{align}
i((-r,q)^T)_t &= L_+^2(r,q)^T =-1/2 \left((1/2)
r_{xx}-r^2q, {(1/2) q_{xx}-q^2r}\right)^T,
 \notag\\
i((-r_k,q_k)^T)_t &=(L+L^{-1})( r_k, q_k)^T= (1-r_kq_k)(E^++E^-)(r_k,q_k)^T. \notag
\end{align}
(Note that the later equation is the differential-difference version of NLS (up to a linear term) arising in the work of Chiu and Ladik, just presented in Eq. \eqref{ALrq} in the previous section. As pointed out earlier, we will also refer to this equation as the AL equation.)

More precisely, setting $\cl_+ \doteq L+L^{-1}$, we note that, based on Eqs. \eqref{L}, \eqref{Linv},

\begin{align}
\cl_+ = (E^++E^-) & +   (1-r_kq_k) \bpm p(\eu-\ed) r_k\\
-\check{p}(\eu-\ed) q_k\epm \notag \\
& +  (1-r_kq_k)(E^+-E^-)\bpm
 {r_k\jkp\left(\frac{q_j}{1-r_jq_j}\right)} &
{-r_k \jkp\left(\frac{r_j}{1-r_jq_j}\right)} \\
{-q_k \jkp\left(\frac{q_j}{1-r_jq_j}\right)} &
 {q_k \jkp\left(\frac{r_j}{1-r_jq_j}\right)}\epm
 \notag \\ & +  (1-r_kq_k)E^-\bpm
 \frac{r_kq_k}{1-r_kq_k} & \frac{-r_k^2}{1-r_kq_k} \\
 \frac{-q_k^2}{1-r_kq_k} & \frac{q_kr_k}{1-r_kq_k}\epm\notag
\\ & +  \bpm r_k \jkp(q_j) & r_k \jkp(r_j) \\
q_k \jkp(q_j) & q_k \jkp(r_j)\epm (E^+-E^-)+
\bpm -r_kq_kE^+ & r_k^2 E^- \\
q_k^2 E^- & -q_kr_k E^+\epm. \label{Lscript+}
\end{align}
\begin{rmk}
From this point on, we will assume that the integration constants $p,\check{p}$ vanish, unless explicitly noted otherwise.
\end{rmk}
\begin{prop}
The operator $\cl_+$ generates a sub-hierarchy of the family of evolution equations \eqref{PotEvol} given by
\begin{equation}
((-r_k,q_k)^T)_t= {\mathcal P}(\cl_+)(r_k,q_k)^T, \label{scripL+Evol}
\end{equation}
where ${\mathcal P}$ denotes an arbitrary polynomial in $\cl_+$.
\end{prop}
\pf  Since $L$ and $L^{-1}$ commute, all powers of $\cl_+=L+L^{-1}$ can also be written as linear combinations of powers of $L$ and $L^{-1}$, and hence ${\mathcal P}(\cl_+)=\om_1(L)+\om_2(L^{-1})$.\epf

Note that the hierarchy \eqref{scripL+Evol} encompasses various versions of the AL equations.  In the reduction $r=-\overline q$, these include \eqref{ALqqbar}, and also the more standard
\begin{equation}
i(q_k)_t=(1+|q_k|^2)(q_{k+1}+q_{k-1})-2q_k,\label{ALstandard}
\end{equation}
obtained from ${\mathcal P}(\cl_+)=\cl_+ -2$.  Observe also the parallel with the Calogero--Degasperis hierarchy \eqref{CalogDegasEvol} associated with NLS.

In the continuous setting, we explored a second recursive-type operator giving rise to the NLS hierarchy generated by $L_+$. This operator, denoted by $R$, was essentially constructed as the sum of $L_+$ and $L_-$ (see \eqref{contL+} and \eqref{contL-}) and was later seen to carry important geometric information concerning the Poisson geometry of NLS.

Our next goal is to construct discrete analogs of the operators $L_-$ and $R$, and to exhibit their connection to the AL hierarchy \eqref{scripL+Evol} stemming from $\cl_+$. The precise geometric
implications of these constructions will be addressed in Section \ref{geom}.

We begin by observing that the form of $L_-$ (specifically the presence of the integral operator $I_-$ in place of the $I_+$ appearing in $L_+$) suggests an adjoint-type relation between $L_+$ and $L_-$.
\begin{prop}\label{PROP:L+L-=adjoints}
Let $\sigma_3\doteq\diag(1,-1)$, and let $^*$ denote the adjoint of a given operator with respect to the inner product \eqref{contInnP}. Then $L_-=\sigma_3L_+^*{\sigma_3}^{-1}.$
\end{prop}
\pf We first compute $L_+^*$ using the aforementioned inner product and integration by parts.  The proposition then follows upon conjugating the calculation's result by $\sigma_3$ \cite{guada}. \epf

In the remaining part of this section we construct the adjoint of $\cl_+$ and use it to define $\cl_-$, the discrete counterpart of $L_-$.  We then show that, once again, the sum $\cl_++\cl_-$ (appropriately conjugated) defines a (discrete) recursion operator $\cR$ for AL which is analogous to the ``continuous" $R$ associated to NLS.

Let $\langle\cdot,\cdot\rangle$ denote the (non-Hermitean) inner product on the complex space
$\cdots\mathbb{C}^2\times\mathbb{C}^2\times\mathbb{C}^2\cdots$ (indexed by $k$ in ${\mathbb Z}$) defined by
\begin{equation}\vspace{-.05in}
\left<(a,b)^T,(c,d)^T\right>
 = -\sum \left(a_k d_k + b_k c_k\right)
(1-q_kr_k)^{-1}, \label{DiscInnProd}
\end{equation}
where $(a, b)^T\doteq(a_k, b_k)^T$, $(c, d)^T\doteq(c_k, d_k)^T$,
$k\in {\mathbb Z}$.

Inspired by the preceding proposition we define
\begin{equation}
\cl_- \doteq \sigma_3\cl_+^*\sigma_3^{-1}=\sigma_3(L^*+(L^{-1})^*)\sigma_3^{-1},
\label{defLscrip-}
\end{equation}
where $^*$ denotes the adjoint of an operator acting on $\cdots\mathbb{C}^2\times\mathbb{C}^2\times\mathbb{C}^2\cdots$ with respect to the inner product \eqref{DiscInnProd} just defined.  We then have,
\begin{align}
L^* & =  (1-r_k q_k)\bpm E^-(\frac{1}{1-r_kq_k}) & 0 \\
0 & E^+(\frac{1}{1-r_kq_k})\epm \nonumber \\ & +  (1-r_kq_k)\bpm
 {E^-\left(r_k \jkm\left(\frac{q_j}{1-r_jq_j}\right)\right)} &
 {E^-\left(r_k \jkm\left(\frac{r_j}{1-r_jq_j}\right)\right)} \\
{-E^+\left(q_k \jkm\left(\frac{q_j}{1-r_jq_j}\right)\right)} &
{-E^+\left(q_k \jkm\left(\frac{r_j}{1-r_jq_j}\right)\right)}\epm
 \nonumber \\ & +  (1-r_kq_k)\bpm
{-E^-\left(\frac{r_kq_k}{1-r_kq_k}\right)} & 0 \\
 {E^+\left(\frac{q_k^2}{1-r_kq_k}\right)} & 0\epm\nonumber
\\ & +  \bpm r_k \jkm(q_jE^-) & r_k \jkm(r_jE^+) \\
-q_k \jkm(q_jE^-) & -q_k \jkm(r_jE^+)\epm +  \bpm
0 & -r_k^2 E^+ \\ 0 & q_kr_k E^+\epm, \label{L*}
\end{align}\vspace{-.09in}
\begin{align}
(L^{-1})^* & =  (1-r_k q_k)\bpm E^+(\frac{1}{1-r_kq_k}) & 0 \\
0 & E^-(\frac{1}{1-r_kq_k})\epm \nonumber \\ & +  (1-r_kq_k)\bpm
{-E^+\left(r_k \jkm\left(\frac{q_j}{1-r_jq_j}\right)\right)} &
{-E^+\left(r_k \jkm\left(\frac{r_j}{1-r_jq_j}\right)\right)} \\
 {E^-\left(q_k \jkm\left(\frac{q_j}{1-r_jq_j}\right)\right)} &
 {E^-\left(q_k \jkm\left(\frac{r_j}{1-r_jq_j}\right)\right)}\epm
 \nonumber \\ & +  (1-r_kq_k)\bpm 0 &
 {E^+\left(\frac{r_k^2}{1-r_kq_k}\right)} \\ 0 &
{-E^-\left(\frac{q_kr_k}{1-r_kq_k}\right)}\epm\nonumber
\\ & +  \bpm -r_k \jkm(q_jE^+) & -r_k \jkm(r_jE^-) \\
q_k \jkm(q_jE^+) & q_k \jkm(r_jE^-)\epm +  \bpm
r_kq_k E^+ & 0 \\ -q_k^2 E^+ & 0\epm,\label{Linv*}
\end{align}
where$\jkm$ denotes the equivalent of the integral operator $I_-$ in the discrete context, that is,$\jkm (u_j)\doteq \sum_{j=-\infty}^k u_j.$ Summing the two expressions above and conjugating by $\sigma_3$,
yields:
\begin{align}
\cl_- = (E^++E^-) & -  (1-r_kq_k)(E^+-E^-)\bpm
 {r_k\jkm\left(\frac{q_j}{1-r_jq_j}\right)} &
{-r_k \jkm\left(\frac{r_j}{1-r_jq_j}\right)} \\
{-q_k \jkm\left(\frac{q_j}{1-r_jq_j}\right)} &
 {q_k \jkm\left(\frac{r_j}{1-r_jq_j}\right)}\epm
 \nonumber \\ & -  (1-r_kq_k)E^+\bpm
 \frac{-r_kq_k}{1-r_kq_k} & \frac{r_k^2}{1-r_kq_k} \\
 \frac{q_k^2}{1-r_kq_k} & \frac{-q_kr_k}{1-r_kq_k}\epm\nonumber
\\ & -  \bpm r_k \jkm(q_j) & r_k \jkm(r_j) \\
q_k \jkm(q_j) & q_k \jkm(r_j)\epm (E^+-E^-) - \bpm r_kq_kE^- & -r_k^2 E^+ \\
-q_k^2 E^+ & q_kr_k E^-\epm. \label{Lscript-}
\end{align}
The (lengthy but straight-forward) calculations yielding formulas \eqref{L*} for $L^*$ and
\eqref{Linv*} for $(L^{-1})^*$ are outlined in Appendix B of \cite{guada}.

Let us now introduce the operator
\begin{align}
\cR\ & \doteq  \bpm 0 & 1 \\ -1 & 0\epm (\cl_++\cl_-)\bpm 0 & -1 \\
1 & 0\epm = 2(E^++E^-) \nonumber \\& +  (1-r_kq_k)(E^+-E^-)\bpm
 {q_k(\jkp-\jkm)\left(\frac{r_j}{1-r_jq_j}\right)} &
 {q_k(\jkp-\jkm)\left(\frac{q_j}{1-r_jq_j}\right)} \\
 {r_k(\jkp-\jkm)\left(\frac{r_j}{1-r_jq_j}\right)} &
 {r_k(\jkp-\jkm)\left(\frac{q_j}{1-r_jq_j}\right)}\epm
 \nonumber \\ & +  (1-r_kq_k)(E^++E^-)\bpm
 \frac{q_kr_k}{1-r_kq_k} & \frac{q_k^2}{1-r_kq_k} \\
 \frac{r_k^2}{1-r_kq_k} & \frac{r_kq_k}{1-r_kq_k}\epm \nonumber
\\ & +  \bpm q_k (\jkp-\jkm)(r_j) & -q_k (\jkp-\jkm)(q_j) \\
-r_k (\jkp-\jkm)(r_j) & r_k (\jkp-\jkm)(q_j)\epm
(E^+-E^-)\nonumber \\ & +  \bpm -q_kr_k & -q_k^2\\
-r_k^2 & -r_kq_k\epm (E^++E^-). \label{discR}
\end{align}
Based solely on its form and method of construction, $\cR$ may be regarded as a discrete analog of $R$.  The analogy is strengthened by the fact that recursive applications of $\cR$ to the
(Hamiltonian) field $(iq_k,-ir_k)^T$ yield the same hierarchy of evolution equations generated by $\cl_+$.  This can be verified by direct calculation for the first few iterations; it will be proven in the general case in Theorem \ref{THM:kc}. For the moment, note --once again-- the parallel with the
continuous setting.

In Section \ref{geom} we will also prove that $\cR$ is in fact the
true \emph{geometric} analog of the recursion operator $R$
associated to the the NLS hierarchy.  In particular, we will
realize $\cR$ as the composition of two skew-symmetric operators
$\cj^{-1}$ and $\ck$ on the appropriate space of complex, vector
valued functions of a discrete real variable $k$.

The next section contains a pivotal result for rigorously proving
the geometric character of $\cR$; thus the framework for unveiling
the Poisson geometry of the AL equations begins to emerge.

\subsection{Resolvent identities for $L$ and $L^{-1}$}

In this section we establish two resolvent-type identities.  The
first one links the hierarchy of evolution equations generated by
powers of $L$ to the hierarchy of gradient fields
\eqref{aGrads-psi2} stemming from the generating function
$\delta_k\log a+\delta_k H_0$.  The second one, relates the flows
given in terms of powers of $L^{-1}$ to the hierarchy of gradient
fields \eqref{ahGrads-psi1} stemming from the generating function
$\delta_k\log \ah+\delta_k H_0$.

The fact that the operators $L$ and $L^{-1}$ satisfy these
identities will be key for postulating and proving several
fundamental results in Section \ref{geom} .  The first of these is
Theorem \ref{THM:kc}, followed by Corollary \ref{COR:conj1} and Corollary
\ref{COR:CDegDisc}.  Through them, we show: (1) that the hierarchy
of flows \eqref{item-iii} generated by the recursion operator
$\cR$  is equivalent to the hierarchy arising from powers of
$\cl_+=L+L^{-1}$; and (2) that the Calogero-Degasperis scheme
leading to Eq. \eqref{L+hier=L-hier} (the continuous ``kernel
condition") and Proposition \ref{PROP:contL+hier=contRhier}
extends to the discrete setting.  Furthermore, once Theorem
\ref{THM:kc} is established, we use the resolvent formulas to
prove that the $\cR$-generated flows are bi-Hamiltonian (Theorem
\ref{THM:biham}.)

In the study of hierarchies of integrable evolution equations, it is morally expected that the operator playing the role of $L$ or $L^{-1}$ (or even $L_+$ in the continuous case,) will satisfy an identity with the properties just described.  This is because the generating function $\delta_k\log a$ (or its equivalent) comprises a full hierarchy of gradient fields which we know should be recursively obtainable from $L$ (or its corresponding counterpart.) In essence, the resolvent identity allows us to compactly write the same hierarchy of gradient fields but now in terms of $L$ and its powers.

Based on this line of reasoning, one expects that the resolvent acting on the simplest, initial flow in the $L$-hierarchy (or $L^{-1}$-hierarchy) should yield the full sequence of gradients. It turns out that, indeed,
\begin{thm}\label{THM:resol}  Let $L$, $L^{-1}$ be as in \eqref{L}, \eqref{Linv} respectively.  Then $L$ and $L^{-1}$ satisfy the following resolvent identities
\begin{align}
(I-z^{-2}L)^{-1}(r_k, q_k)^T&=(1-r_kq_k)(\delta_k \log a+\delta_k H_0), \label{resolvL} \\
(I-z^2L^{-1})^{-1}(r_k,q_k)^T&=(1-r_kq_k)(\delta_k \log \ah+\delta_k H_0), \label{resolvLinv}
\end{align}
where $H_0=-\sum_{k=-\infty}^\infty \log(1-r_kq_k).$
\end{thm}
\pf  Setting $(\bbk,-\aak)^T=P(k+1)/a\ (\sit_k\fit_{k+1},-\fio_k\sio_{k+1})^T$,
in accordance with \eqref{WkMat} in Appendix II, we recall from Theorem \ref{THM:grads} that the non-standard squared eigenfunctions given by $-\alpha_k$ and $\beta_k$ encode the hierarchy of gradient fields stemming from the generating function  $\ddk \log a + \ddk H_0$. The resolvent formula \eqref{resolvL} we seek to establish explicitly portrays $(-\aak,  \bbk)^T$ as a generating function for the hierarchy of fields arising from powers of $L$.

Guided by the form of the leading term of $L$  we select the first component of \eqref{-sh} and the second component of \eqref{+sh} (re-writing the $\delta\gamma$-factor in each using \eqref{abgd2}) and combine them so as to get
\begin{multline}
\bpm \ed & 0\\0 & \eu\epm \left[(1-r_kq_k) \abv \right]= z^2\left[(1-r_kq_k)\abv-\bpm r_k\\q_k\epm\right]\label{+-sh} \\
+\bpm -(r_{k-1}+z^2r_k) &0 \\ 0& {q_{k+1}+z^2q_k}\epm\ \bpm \ed & 0\\0 &1\epm \bpm  -J^+_{k+1}(q_j) & J^+_{k+1}(r_j)\\ J^+_{k+1}(q_j) & -J^+_{k+1}(r_j) \epm\abv.
\end{multline}
Re-grouping terms and dividing through by $-z^2$, the last equation is equivalent to
\begin{equation}\left(I-z^{-2}\bpm \ed & 0 \\ 0 & \eu\epm
\right)\left[(1-r_kq_k)\abv\right]+z^{-2}{{\mathcal B}}=\bpm r_k \\ q_k\epm, \label{idea}
\end{equation}
where ${{\mathcal B}}$ is the second term on the right-hand-side of \eqref{+-sh}.

Based on the expected form of the resolvent identity and the formula for $L$, Eq. \eqref{idea} suggests that ${{\mathcal B}}$ ought to simply amount to all but the first term of the operator $-L$ acting on the vector $[(1-r_kq_k)(\bbk, -\aak)^T]$.  In other words,
\begin{align}
{{\mathcal B}}&\?\bpm -r_k & 0 \\ 0 & q_k\epm \Upsilon(k)\bpm \ed & 0 \\ 0
&\eu\epm \left[(1-r_kq_k)\abv\right]\nonumber \\
&+ \bpm -r_{k-1} & 0 \\ 0 & q_{k+1}\epm
(1-r_kq_k)\Upsilon(k)\abv, \ \ \ \Upsilon(k)\doteq\bpm -\jkp(q_j) & \jkp(r_j) \\
                J^+_{k+1}(q_j) & -J^+_{k+1}(r_j) \epm,\label{claim&upsilon}
\end{align}
(where the symbol $\?$ signifies a claim yet-to-be-established.) Expressing ${{\mathcal B}}$ (in \eqref{+-sh}) in terms of $\Upsilon(k)$, one notices that \eqref{claim&upsilon} is equivalent to $\diag(-r_k, q_k){{\mathcal A}}\?0$, where
\begin{equation}
{{\mathcal A}}= \left[\bpm r_{k-1}q_k+z^2 & 0 \\ 0 &
q_{k+1}r_k+z^2\epm\Upsilon(k)-\Upsilon(k)\bpm \ed & 0 \\ 0 & \eu
\epm(1-r_kq_k)\right]\abv,\label{Adiff}
\end{equation}
which in turn holds {\it if} ${{\mathcal A}}=0$.  Iterating once, we re-write the factor
$\diag(\ed,\eu)[(1-r_kq_k)(\bbk,-\aak)^T]$ in the second term of ${{\mathcal A}}$ with its equivalent form given in Eq. \eqref{+-sh}.   Grouping terms in $z^2$, we see that ${{\mathcal A}}$ amounts to the following quadratic equation in $\Upsilon(k)$:
\begin{align}
&\left[\bpm r_{k-1}q_k & 0 \\ 0 &
q_{k+1}r_k\epm- \Upsilon(k) \bpm -r_{j-1}& 0 \\
0 & q_{j+1}\epm\right]\left[\Upsilon(k)\abv\right] \notag \\
& + z^2 \Upsilon(k)\left[r_jq_j\abvj-\bpm -r_j & 0 \\ 0 & q_j \epm
\bpm
q_j\beta_j+r_j\alpha_j \\ 0\epm\right].\label{ah!}
\end{align}
Now, due to the form of the matrix $\Upsilon(k)$, a simple calculation shows that the coefficient of the $z^2$ term in Eq. \eqref{ah!} vanishes.  The quadratic term within the first line of \eqref{ah!} telescopes down to its boundary term, $\jkp (q_j\beta_j+r_j\alpha_j)(q_kr_{k-1}, -q_{k+1}r_k)^T$, which is evidently the negative of the linear term in the first line of \eqref{ah!}.  Hence ${{\mathcal A}}=0$ and the resolvent formula follows. \epf

\section{Geometry: almost Poisson structure of AL}\label{geom}

In this section we exhibit the bi-Hamiltonian character of the AL hierarchy of equations, as defined in Eq.
\eqref{item-iii}.

This postulated characterization is attained through the construction of two geometrically meaningful operators on the phase space $\cc$ of these equations: an almost-Poisson operator $\ck$, and a recursion operator $\cR$ which can be both realized in terms of $\cl_+$ and $\cl_-$ (see \eqref{Lscript+} and \eqref{Lscript-}), and as $\ck\cj^{-1}$ (where $\cj=\diag(-i,i)$ is the standard Poisson operator on $\cc$.)

The former depiction of $\cR$, allows us to view the $\cR$-generated hierarchy \eqref{item-iii} as a sub-hierarchy of the Chiu--Ladik flows.  The latter realization of this recursion operator, attests to its true geometric character, and it is responsible for elucidating the bi-Hamiltonian character of the hierarchy.

The skew operator $\ck$ can be portrayed as the discrete counterpart of Magri's Poisson operator for NLS, based both on its mode of construction and on its geometric properties.  For example, we will see that, once their bi-Hamiltonian nature is established, the skewness of $\ck$ serves to prove the commutativity of the AL flows generated by $\cR$.

\subsection{The Geometric Context: AL as a Hamiltonian System}

We begin by defining the phase space $\cc$ of the AL equations as the space of complex, vector valued functions of a discrete real variable $k\in{\mathbb Z}$ given by
$$
\cc=\{(q_k,r_k)^T \in \cdots\times{\mathbb C}^2\times{\mathbb C}^2\times{\mathbb C}^2\times\cdots :
|q_k|,|r_k|\rightarrow 0 \mbox{ as } k\rightarrow\pm\infty, \}.
$$
A point in $\cc$ will either be denoted by $ (q,r)^T\doteq \left[ \cdots, (q_k,r_k)^T, (q_{k+1},r_{k+1})^T, \cdots \right],$ or just by $(q_k, r_k)^T$, depending on the context.

As in the continuous case, we may endow the tangent-bundle to $\cc$, with a non-degenerate bilinear form $\langle \cdot,\cdot \rangle$ arising from a point-dependent inner product defined by
\begin{equation}
\left<(a,b)^T,(c,d)^T\right>_{(q, r)^T}= -\sum \left(a_k d_k + b_k c_k\right) (1-q_kr_k)^{-1}, \label{DiscRieMetric}
\end{equation}
on the tangent space to $\cc$ at $(q, r)^T$; that is, $(a, b)^T, (c, d)^T$ lie in $T_{(q, r)^T}\cc\simeq \cc$.  Under the identification $r_k=-\overline{q}_k$, \eqref{DiscRieMetric} defines a positive-definite, real inner product.

Given a functional $F$ over $\cc$ we may define the (discrete) functional gradient of $F$, $\nabla F$, by
$\nabla F ((q,r)^T)=[ \cdots, \nabla_k F((q,r)^T), \nabla_{k+1} F((q,r)^T), \cdots],$ where for $k$ in ${\mathbb Z}$,
\begin{equation}
\nabla_k F((q,r)^T) =  -(1-q_kr_k)\sigma_1(\delta F/\delta q_k, \delta F/\delta r_k)^T \doteq \nabla F ((q_k, r_k)^T),\label{discGrad}
\end{equation}
where $\sigma_1$ is the usual Pauli matrix with {\it ones} in the off-diagonal.

Now, let $\cj_{(q, r)^T}\doteq\diag(\cdots,\cj_k, \cj_{k+1}, \cdots)$, $\cj_k=\diag(-i,i)$ for all $k$,
define the (standard, point-independent) skew-symmetric operator on $T_{(q, r)^T}\cc$, acting as multiplication by $\cj_k$ on each ${\mathbb C}^2$ component.  The operator $\cj$ (which may also be
regarded as an operator on the tangent bundle to $\cc$) may be used to define the following Poisson bracket on $\cc$:
\begin{align}
&\{F, G \}_\cj((q, r)^T) = \langle \nabla F, \cj \nabla G\rangle((q, r)^T)\notag \\
&= \left< (1-q_kr_k)(\delta F/\delta r_k,\delta F/\delta q_k)^T, \diag(-i,i)(1-q_kr_k)(\delta F/\delta r_k,\delta F/\delta q_k)^T \right>_{(q, r)^T},\label{discJbracket}
\end{align}
where, $F$ and $G$ are smooth functionals over $\cc$ which become real-valued under the identification $r_k=-\overline{q}_k$.  Note that as in the continuous setting, this reality requirement
implies that $F$ and $G$ must be either symmetric or anti-symmetric in $r_k$, $q_k$ (for all $k$.)
\begin{rmk}
Since the operator $\cj_k$ defined above does not actually depend on $k$, we will often use the notation $\cj$ to refer to the $2\times 2$ complex matrix $\mbox{Diag}(-i, i)$ originally defined as $\cj_k$.  By the same token, from now own we may omit noting explicitly the point dependence of operators which are clearly characterized as such by their definitions.
\end{rmk}
Let $H((q_k,r_k)^T)=-\sum_{-\infty}^\infty (q_kr_{k+1}+r_kq_{k+1})$.  Then the AL equations \eqref{ALrq} can be obtained as a Hamiltonian system on $\cc$ with respect the $\cj$-bracket \eqref{discJbracket} for the Hamiltonian functional $H$.  In this same manner, the Hamiltonian $ H((q_k,r_k)^T)=-\sum_{-\infty}^\infty (q_kr_{k+1}+r_kq_{k+1}+2\log(1-q_kr_k))$ produces the standard (vector) AL flow, which becomes \eqref{ALstandard} in the reduction $r_k=-\overline{q}_k$. (See also \cite{FT}.)

\subsection{The Operators $\ck$ and $\cR$ and the AL Hierarchy}

Guided by the parallel with our continuous setup and the form of the operators $\cl_+$ and $\cl_-$, we re-consider the operator $\cR$ defined in \eqref{discR}, and define a new operator
\begin{align}
\ck ((q_k,r_k)^T)
&\doteq \cR((q,r)^T)\cj = -i\Bigg\{ 2(E^++ E^-) \notag \\
&+ (1-r_k q_k)(E^+-E^-)\bpm              q_k (J_k^+ - J_k^-)
        \left( \frac{r_j}{1-r_jq_j}\right)&    -q_k (J_k^+ - J_k^-)
        \left( \frac{q_j}{1-r_jq_j}\right)\\    r_k (J_k^+ - J_k^-)
        \left( \frac{r_j}{1-r_jq_j}\right)&   -r_k (J_k^+ - J_k^-)
        \left( \frac{q_j}{1-r_jq_j}\right)\epm \notag\\
&+ (1-r_k q_k)(E^++E^-)
        \bpm  \frac{q_k r_k}{1-r_k q_k}&  \frac{-q_k^2}{1-r_k q_k}
        \\    \frac{r_k^2}{1-r_k q_k}&  \frac{-r_k q_k}{1-r_k q_k}\epm \notag\\
&+ \bpm                q_k (J_k^+ - J_k^-)
        \left(r_j\right)&   q_k (J_k^+ - J_k^-)
        \left(q_j\right)\\  -r_k (J_k^+ - J_k^-)
        \left(r_j\right)&    -r_k (J_k^+ - J_k^-)
        \left(q_j\right)\epm (E^+-E^-)\notag\\
&+ \bpm -q_k r_k & q_k^2\\ -r_k^2 & r_k q_k\epm (E^++E^-)
        \Bigg\}. \label{discK}
\end{align}
Just as with $\cR$, $\ck$ is point-dependent and non-local; its
definition at site $k$ involves functions of the vector potential
at neighboring sites, for example, functions of $q_{k\pm1}$,
$r_{k\pm1}$.

As $\ck$ acts on the tangent bundle of $\cc$, it may be used to define the bracket
\begin{equation}
\{F, G \}_\ck((q_k, r_k)^T)= \langle \nabla F, \ck \nabla G\rangle((q,r)^T),\label{discKbracket}
\end{equation}
on the class of functionals over $\cc$ just discussed.
\begin{thm}
The operator $\ck$ is skew-symmetric.
\end{thm}
\pf  One first shows that $\cR^*=\sigma_3\cR\sigma_3$
\cite{guada}, where, as before, $^*$ denotes the adjoint of an
operator with respect to the inner product \eqref{DiscRieMetric},
and $\sigma_3=\mbox{Diag}(1, -1)$.  This calculation uses the fact
that both $\pm i\sigma_2)$ and $\sigma_3$ are self-adjoint, while
$\cl_-\doteq \sigma_3 \cl_+^*\sigma_3^{-1}$.

The result now follows from the definition of $\ck$ in terms of $\cR$ and $\cj=-i\sigma_3$.  Indeed,
\begin{align}
\{F,G\}_\ck & =  \langle \nabla F,\ck \nabla G\rangle = \langle
                  \nabla F,\cR\cj \nabla G\rangle =\langle \cR^* \nabla F,\cj \nabla G\rangle= \langle -\cj
                  (\sigma_3\cR\sigma_3) \nabla F, \nabla G\rangle \notag \\
            & =  \langle\cR(i\sigma_3) \nabla F, \nabla G\rangle =  -\langle\cR\cj \nabla F, \nabla G\rangle =
                  -\langle \nabla G, \ck \nabla F \rangle = -\{G,F\}_\ck.  \ \ \ \ \ \ \ \ \ \ \ \ \ \square \notag
\end{align}
\begin{cor} The $\ck$-bracket \eqref{discKbracket} defines an almost Poisson structure on $\cc$ .\end{cor}
\pf  In addition to being skew, the $\ck-bracket$ satisfies the Leibnitz identity, as we show in \cite{guada}.
\epf

\subsection{The bi-Hamiltonian character of the AL equations}

Let us now consider the functional $H_0((q_k,r_k)^T)=-\sum_{-\infty}^\infty \log(1-q_kr_k),$ on $\cc$, and observe that $\cj\nabla H_0((q_k, r_k)^T)=i(q_k,-r_k)^T.$  If we now define $X_{(n)}\doteq\cR^n\cj\nabla H_0,$ then $X_{(n)}$ defines a hierarchy of fields on the tangent bundle
to $\cc$ which, we claim, has the following properties:
\begin{itemize}
\item[(\emph{i})] The fields $X_{(n)}$ mutually commute;

\item[(\emph{ii})] The hierarchy comprises the AL equations \eqref{ALrq}, which occur as $X_{(1)}$;

\item[(\emph{iii})] The evolution equations
\begin{equation}
((q_k,r_k)^T)_t = \cR^n\cj\nabla_k H_0=X_{(n)}, \label{item-iii}
\end{equation}
are equivalent to a sub-hierarchy of the family \eqref{scripL+Evol} and are therefore integrable;

\item[(\emph{iv})] The fields $X_{(n)}$ in the hierarchy are Hamiltonian with respect to $\cj$ and also with respect to $\ck$, i.e., bi-Hamiltonian.
\end{itemize}

\begin{rmk}\label{RMK:Lenard}
As indicated in \cite{guada}, the term \emph{bi-Hamiltonian} here
refers to sequences of fields $X_{(n)}$ which satisfies the
\emph{Lenard relations} $X_{(n)}  =\ck\nabla H_{n-1} =\cj\nabla
H_n$ \cite{Damianou}.  The clarification is pertinent due to the
fact that the Poisson nature of the skew operator $\ck$ has not
yet been established.
\end{rmk}

It is a good exercise to verify items (\emph{iii}) and (\emph{iv}) explicitly for the first 2 iterations \cite{guada} (in the process, item (\emph{ii}) becomes apparent.)  As such calculations indicate,
being able to show that the fields $\cl_+^n((r_k,q_k)^T)$ are in the kernel of $\cl_+-\cl_-$ for all $n$ is of central importance in proving claim (\emph{iii}).  This is precisely the content of Theorem \ref{THM:kc}.  In it, we use once again the identities derived in Appendix II to extend the Calogero-Degasperis ``kernel-condition" result
(see \eqref{L+hier=L-hier}) to the discrete setting.  Corollary \ref{COR:conj1}, the discrete-analog of Proposition \ref{PROP:contL+hier=contRhier}, will establish item (\emph{iii}) in full generality.  Theorem \ref{THM:biham} will then use the resolvent identities \eqref{resolvL}, \eqref{resolvLinv} together with Corollary \ref{COR:conj1} to establish (\emph{iv}). The proof of item (\emph{i}) is rather immediate given (\emph{iv}).  It is the content of Corollary \ref{COR:cflows}.

Throughout the remainder of this note, $\sum f_j=\sum_{k=-\infty}^{\infty}$; also the difference operator $\diff$ will often be written as $\diff=\Do+\Dt,$ where
\begin{equation}
\Do\doteq (1-r_kq_k)(E^+-E^-)
                  \bpm r_k\sum
                  q_j/(1-r_jq_j)&
                  -r_k\sum
                  r_j/(1-r_jq_j) \\  -q_k\sum
                  q_j/(1-r_jq_j) &
                  q_k\sum
                  r_j/(1-r_jq_j)\epm,
\label{D1}
\end{equation}
and
\begin{equation}
\Dt\doteq \bpm r_k\sum q_j &
                  r_k\sum r_j \\  q_k\sum
                  q_j& q_k\sum r_j\epm
                  (E^+-E^-).
\label{D2}
\end{equation}
\begin{lem}\label{lem1}
Let $j=1,2$.  If  $\D_j\left(L^m(r_k,q_k)^T\right)=\D_j\left({(L^{-1})}^m(r_k,q_k)^T\right)=0$ for all $m$, then $(\diff)\left(\cl_+^n (r_k,q_k)^T\right)=0$ for all $n$ ($m,n$ non-negative integers.)
\end{lem}
\pf    This follows immediately from the definition of $\cl_+$ (see \eqref{Lscript+} above,) the inverse relation between $L$ and $L^{-1}$ (described in \cite{guada}, proven in \cite{LC}) and the binomial theorem.
\epf
\begin{prop}\label{prop1}
For $\Do$ as in \eqref{D1}, $\Do\left(L^m(r_k,q_k)^T\right)=\Do\left({(L^{-1})}^m(r_k,q_k)^T\right)=0.$
\end{prop}
\pf  Recall that the vector of squared eigenfunctions $\abvr$ defined in Appendix II may be explicitly viewed a generating function for the hierarchy of fields arising from powers of $L$ (through the resolvent identities in Theorem \ref{THM:resol}.) The same is true of $\abvhr$ and the fields given by $L^{-1}((r_k, q_k)^T)$.  Using the inner product \eqref{DiscRieMetric} one observes certain relations between the aforementioned fields, which can be compactly expressed as:
\begin{align}
\langle (1-r_kq_k)\abvr, -i\sigma_3(1-r_k q_k)\abvr \rangle &=  i \sum (1-r_kq_k)(\beta_k\alpha_k-\alpha_k\beta_k) = 0, \label{x}\\
\langle (1-r_kq_k)\abvhr, -i\sigma_3(1-r_k q_k)\abvhr\rangle &= 0. \label{xx}
\end{align}
We note that the previous identities translate into commutativity relations among special sub-classes of Hamiltonian functions, namely $\{H_m^a, H_n^a\}=\{H_m^{\ah}, H_n^{\ah}\}=0$ (refer to \cite{guada} for the precise definitions of $H_n^a$ and $H_n^{\ah}$.)  In addition, one may focus attention on the particular relations
\begin{align}
0&=\langle(r_k,q_k)^T,-i\sigma_3(1-r_kq_k)\abvr \rangle = i\sum(q_k\beta_k+r_k\alpha_k), \label{first}\\
0&=\langle (r_k,q_k)^T,-i\sigma_3(1-r_kq_k)\abvhr \rangle =-i\sum(q_k\betah_k+r_k\alphah_k),  \label{second}
\end{align}
obtained from expanding the first entry on the left-hand-side of \eqref{x} (respectively \eqref{xx}) in powers of $z^{-2}$ (respectively $z^2$) as suggested by the resolvent formula \eqref{resolvL} (respectively \eqref{resolvLinv}.)

Equation \eqref{first} and \eqref{D1} together with \eqref{second} and \eqref{D1} immediately show that
\begin{equation}
\Do\left( (1-r_kq_k)\abvr\right)=0, \ \ \mbox{ and } \ \ \Do\left( (1-r_kq_k)\abvhr\right)=0. \label{third/fourth}
\end{equation}
Using now \eqref{resolvL} and the first part of \eqref{third/fourth}, we see that
$$
0=\Do(r_k,q_k)^T+z^{-2}\Do\left( L(r_k,q_k)^T\right)+\cdots+z^{-2m}\Do\left( L^m(r_k,q_k)^T\right)+\cdots
$$
or, equivalently, $\Do\left(L^m(r_k,q_k)^T\right)=0$, for all $m$. The identity $\Do\left( (L^{-1})^m(r_k,q_k)^T\right)=0$ stems from Eqs. \eqref{resolvLinv} and the second piece of \eqref{third/fourth} in  a parallel manner.
\epf
\begin{thm}\label{THM:kc}
For any non-negative integer $n$, $(\diff)\left(\cl_+^n(r_k,q_K)^T\right)=0.$
\end{thm}
\pf  Due to Lemma \ref{lem1} and Proposition \ref{prop1} it suffices to prove that
\begin{equation}
\Dt\big({(L^{-1})}^m(r_k,q_k)\big)\ida\Dt\big(L^m(r_k,q_k)\big)\idb0.\label{remain}
\end{equation}
Below, we present the argument showing that equality (b) in \eqref{remain}
holds.   An analogous argument (based on the identities given at the end of Appendix II) shows that identity (a) holds as well, validating the theorem.

The key step in the argument closely resembles the recurrence strategy employed in the derivation of the resolvent identities \eqref{resolvL}, \eqref{resolvLinv}, as described in Appendix II.  In the current scenario, the technique entails re-writing  the vector $(\eu-\ed) \left[ (1-r_kq_k)\abvr\right]$ as
\begin{align}
&(z^{-2}-z^2)\diag(1,-1)\left[(1-r_kq_k)\abvr - (r_k,q_k)^T \right] \label{exp1} \\
+ & \diag({r_{k+1}+z^{-2}r_k},\ {q_{k+1}+z^2q_k})\ ({q_k\beta_k+z\gamma_k+1},\ {-(r_k\alpha_k+z^{-1}\delta_k-1)})^T \label{exp2}\\
- & \diag( -({r_{k-1}+z^2r_k)}, {q_{k-1}+z^{-2}q_k})\  \ed\ ({q_k\beta_k+z\gamma_k}, {r_k\alpha_k+z^{-1}\delta_k})^T,  \label{exp3}
\end{align}
by considering the difference between equations \eqref{+sh} and \eqref{-sh} in Appendix II.

The remainder of the proof amounts to a series of calculations.
First we utilize Eqs. \eqref{exp1}--\eqref{exp3} together with
\eqref{3} and \eqref{abgd2} in Appendix II, to obtain a simplified
three-term expression for $\Dt((1-r_kq_k)(\beta_k, -\alpha_k)^T)$.
We then exploit the particular form of the resulting expressions
to show that their combined sum vanishes.  Finally, the same
argument used at the end of the proof of Proposition \ref{prop1},
shows that the vanishing of the vector  $\Dt((1-r_kq_k)(\beta_k,
-\alpha_k)^T)$ implies the vanishing of $\Dt (L^m (r_k, q_k)^T)$
for all $m$.

Due to \eqref{first}, the matrix pre-factor of the operator $\Dt$ applied to expression \eqref{exp1}, becomes:
\begin{equation}
 (z^{-2}-z^2)\bpm r_k & 0 \\ 0& q_k\epm \bpm \sum q_j & -\sum r_j\\ \sum q_j & -\sum r_j\epm \left[-r_kq_k\abv\right], \label{exp1new}
\end{equation}
Expression \eqref{exp1new} is the first term in the sought-after three-term expression for $\Dt((1-r_kq_k)(\beta_k, -\alpha_k)^T)$.  The second and third terms (appearing in \eqref{exp2new}, \eqref{exp3new}) are essentially obtained by re-writing \eqref{exp2} and \eqref{exp3} by means of \eqref{abgd2}, in Appendix II.  Specifically, the matrix pre-factor of $\Dt$ applied to expression $\eqref{exp2}$, yields the new version of the second term in the expression for $\Dt((1-r_kq_k)(\beta_k, -\alpha_k)^T)$, namely,
\begin{align}
& \bpm r_k & 0\\ 0& q_k\epm\bpm\sum q_j & \sum r_j \\ \sum q_j & \sum r_j \epm \bpm {r_{k+1}+z^{-2}r_k} & 0 \\ 0 & {q_{k+1}+z^2q_k}\epm \notag\\
&\hspace{+.6in}\Bigg[ \bpm  -J^+_{k+1}(q_j) & J^+_{k+1}(r_j) \\ J^+_{k+1}(q_j) & -J^+_{k+1}(r_j) \epm \abv +\bpm 1 \\ 0\epm \Bigg]  \notag \\
= &\diag(r_k,q_k) \Big[ \sum q_j(r_{j+1}+z^{-2}r_j) (1,1)^T \notag \\
 & + {\sum\left[-q_j(r_{j+1}+z^{-2}r_j)+r_j(q_{j+1}+z^2q_j)\right]J^+_{j+1}(q_l\beta_l+r_l\alpha_l)} (1,1)^T
 \Big] \notag \\
= &\diag(r_k,q_k) \sum \big[ q_j(r_{j+1}+z^{-2}r_j) \notag \\
 & + {\left(-q_{j-1}(r_j+z^{-2}r_{j-1})+r_{j-1}(q_j+z^2q_{j-1})\right)J^+_j(q_l\beta_l+r_l\alpha_l)} \big] (1,1)^T \notag \\
= &\diag(r_k,q_k) \sum\big[q_j(r_{j+1}+z^{-2}r_j) \notag \\
& + (q_j\beta_j+r_j\alpha_j) J^-_j\left(-q_{l-1}r_l+r_{l-1}q_l-(z^{-2}-z^2)q_{l-1}r_{l-1}\right)\big] (1,1)^T,  \label{exp2new}
\end{align}
where, recall, $\jkm (u_j)\doteq \sum_{j=-\infty}^k u_j.$ The last equality in the above string of identities results from a direct application of the summation by parts formula (see Appendix A in \cite{guada}.)  The previous one, amounts simply to a convenient shift in the summation index for the bi-infinite sum.  This same series of calculations leads to the desired third and final term.  That is, using \eqref {abgd2}, and summation by parts, the action of the matrix pre-factor of $\Dt$ on \eqref{exp3} leads to the vector
\begin{align}
&\diag(r_k, q_k) \sum \big[ -r_j(q_{j-1}+z^{-2}q_j) \notag \\
& + (q_j\beta_j+r_j\alpha_j) J^-_j\left(-r_{l-1}q_l+q_{l-1}r_l+(z^{-2}-z^2)q_lr_l\right)\big] (1,1)^T.  \label{exp3new}
\end{align}

Consider now the sum of expressions \eqref{exp2new} and \eqref{exp3new}.  The $J^-_j$-independent terms all disappear: the two $z$-dependent ones cancel out; the remaining two telescope away, leaving no boundary terms.  Focusing now on the $J^-_j$-dependent terms, one observes that the $z$-independent ones cancel off identically.  The four $z$-dependent terms are pairwise-telescoping, and hence reduce to $\ \diag(r_k,q_k)  \sum \big[ (q_j\beta_j+r_j\alpha_j) q_j r_j (z^{-2}-z^2)\big] (1,1)^T, $
which is exactly the negative of \eqref{exp1new}.  It follows that, for all (non-negative) $m$,
\begin{equation}
\Dt\left((1-r_kq_k)\abvr\right)=0, \ \mbox{ and hence, } \ \Dt\left(L^m (r_k,q_k)^T\right)=0.
\end{equation}
A parallel argument, based on identities \eqref{+shHat}, \eqref{-shHat}, \eqref{abgd2Hat} in Appendix II yields
\begin{equation}
\Dt\big((1-r_kq_k)\abvhr\big)=0, \ \mbox{ and hence, } \ \Dt\big({L^{-1}}^m (r_k,q_k)^T\big)=0.
\end{equation}
The theorem follows.
\epf

\begin{cor}\label{COR:conj1}
 If $(\cl_+-\cl_-)\left(\cl_+^m(r_k,q_K)^T\right)=0 \mbox{ for } m\geq 0, \label{kc}$ then
\begin{equation}
\bpm 0 & 1\\-1 & 0\epm(\cl_++\cl_-)^n\bpm 0 & -1\\1 & 0\epm \cj\nabla_k H_0=2^n\bpm 0 & 1\\-1 & 0\epm \cl_+^n\bpm 0 & -1\\1 & 0\epm \cj\nabla_k H_0, \label{c1}
\end{equation}
for $n\geq 0$ and $H_0=-\srl \log(1-r_kq_k)$.
\end{cor}
\pf  The statement may be verified by direct computation for $n=1$.  One may then proceed by induction on $n$.   Assuming that identity \eqref{c1} holds for $n-1$, we show below that it also holds for $n$.  Theorem \ref{THM:kc} is used in the last line below.  Indeed,
\begin{align}
&\bpm 0 & 1\\-1 & 0\epm(\cl_++\cl_-)^n\bpm 0 & -1\\1 & 0\epm \cj\nabla_k H_0=i\bpm 0 & 1\\-1 & 0\epm(\cl_++\cl_-)(\cl_++\cl_-)^{n-1}\rqv\notag \\
&=\bpm 0 & 1\\-1 & 0\epm(\cl_++\cl_-)\Bigg[i\bpm 0 & -1 \\ 1 & 0\epm\bpm 0 & 1 \\ -1 & 0\epm (\cl_++\cl_-)^{n-1}\rqv\Bigg]  \notag \\
&=\bpm 0 & 1\\-1 & 0\epm(\cl_++\cl_-)\bpm 0 & -1 \\ 1 & 0\epm\Bigg[i2^{n-1}\bpm 0 &1 \\ -1 &0\epm \cl_+^{n-1}\rqv\Bigg]  \notag \\
&=i2^{n-1}\bpm 0 & 1\\-1 & 0\epm(\cl_++\cl_-)\cl_+^{n-1}\rqv=i2^{n-1}\bpm 0 & 1\\-1 & 0\epm\Bigg[ \cl_+^n\rqv+\cl_-\left(\cl_+^{n-1}\rqv\right)\Bigg] \notag \\
&=i2^{n-1}\bpm 0 & 1\\-1 & 0\epm\Bigg[ \cl_+^n\rqv+\cl_+^n\rqv\Bigg] =2^n\bpm 0 & 1 \\ -1 &0\epm\cl_+^n\Bigg[\bpm 0&-1\\1&0\epm\cj\nabla_kH_0\Bigg] \notag. \ \ \ \ \ \square
\end{align}

\begin{cor}\label{COR:CDegDisc}
(Validity of Item (iii).) The evolution equations \eqref{item-iii} are equivalent to a sub-hierarchy of the Chiu-Ladik family \eqref{scripL+Evol} and, thus, integrable.
\end{cor}
\pf
Working from Eq. \eqref{item-iii} and using Corollary \ref{COR:conj1}, we have:
\begin{equation}
\bpm q_k \\ r_k\epm_t = \cR^n\cj\nabla_kH_0=i\bpm 0 & 1 \\ -1 &
0\epm (\cl_++\cl_-)^n\bpm 0 & -1 \\ 1 & 0\epm\bpm q_k \\ -r_k\epm= i2^n \bpm 0 & 1 \\ -1 & 0\epm
\cl_+^n\bpm r_k \\ q_k\epm; \nonumber
\end{equation}
or, equivalently, $((-r_k,q_k)^T)_t=i2^n  \cl_+^n(r_k,q_k)^T,$ which is of form \eqref{scripL+Evol}, as claimed.
\epf

\begin{thm}\label{THM:biham}
(Validity of Item (iv).) The fields $X_{(n)}$ in the hierarchy \eqref{item-iii} are Hamiltonian with respect to both, $\cj$ and $\ck$, i.e., bi-Hamiltonian in the sense of Remark \ref{RMK:Lenard}.
\end{thm}
\pf
Let us once again consider the resolvent formulas for $L$ and $L^{-1}$.  Starting from \eqref{resolvL} and matching inverse powers  of $z^2$, we see that:
\begin{equation}
(1-r_kq_k)\delta_k H_0=(r_k, q_k)^T,\ \mbox{ and also, } \ (1-r_kq_k)\delta_k H_j^a=L^j(r_k,q_k)^T. \label{no1}
\end{equation}
Using the definition of $\nabla_k$ given in \eqref{discGrad}, we observe
\begin{align}
\cj\nabla_k H_j^a &= \diag(-i,i)\left[-\bpm 0 & 1\\1 & 0\epm L^j \bpm 0 & 1\\1 & 0\epm
\bpm q_k\\r_k\epm  \right] \nonumber \\
&=i \bpm 0 & 1\\-1 & 0\epm L^j \bpm 0 & -1\\1 & 0\epm\bpm
q_k\\-r_k\epm= \bpm 0 & 1\\-1 & 0\epm L^j \bpm 0 & -1\\1 & 0\epm
\cj\nabla_kH_0. \label{no3}
\end{align}
Reasoning in a parallel fashion for the $L^{-1}$ resolvent, \eqref{resolvLinv}, we get:
\begin{equation}\small{
(1-r_kq_k)\delta_k H_j^{\ah}=(L^{-1})^j \bpm r_k \\ q_k\epm, \ \ \   \cj\nabla_k H_j^{\ah}=\bpm 0 & 1\\-1 & 0\epm (L^{-1})^j \bpm 0 & -1\\1 & 0\epm \cj\nabla_kH_0. \label{no4}
}\end{equation}
(Notice that, based on Eqs. \eqref{oooo} and \eqref{ooooo}, one has $\delta_k H_j^a=\delta_k C_j$, $\delta_k H_j^{\ah}=\delta_k\hat{C}_j$, for $j>0$.)
Based in Eqs. \eqref{no3} and \eqref{no4}, we have:
\begin{align}
\bpm q_k \\ r_k\epm_t &=\cR^n\cj\nabla_kH_0 \stackrel{\tiny{\mbox{Cor.}\ref{COR:conj1}}}{=}2^n\left[\bpm 0 & 1\\-1 & 0\epm (\cl_+)^n \bpm 0 & -1\\1 & 0\epm \cj\nabla_kH_0\right]\notag\\
&=2^n\left[ \bpm 0 & 1\\-1 & 0\epm (L+L^{-1})^n \bpm 0 & -1\\1 & 0\epm \cj\nabla_kH_0  \right]\notag\\
&=2^n \bpm 0 & 1\\-1 & 0\epm \left[\sum_{n-k\geq k}\bpm n\\k\epm L^{(n-k)-k} + \sum_{n-k< k}\bpm n\\k\epm (L^{-1})^{k-(n-k)} \right]\bpm 0 & -1\\1 & 0\epm \cj\nabla_kH_0\notag \\ &= 2^n\cj\nabla_k\left( \sum_{n-k\geq k}\bpm n\\k\epm H_{n-2k}^a + \sum_{n-k< k}\bpm n\\k\epm H_{-(n-2k)}^{\ah}\right)=\cj\nabla_k \tilde{H}_n, \notag
\end{align}
and so, the last equality, explicitly depicts the alluded hierarchy of flows as $\cj$-Hamiltonian.  Given this, one uses the fact that $\cR=\ck\cj^{-1}$ in a recursive fashion to show that these flows are also $\ck$-Hamiltonian (in the sense of Remark \ref{RMK:Lenard}.)  Indeed, for $n=1$ and $n=2$,
$$
\cR\cj\nabla_kH_0=\ck\nabla_kH_0=\cj\nabla_k\tilde{H}_1, \ \ \mbox{ and } \ \ \cR^2\cj\nabla_kH_0=\ck\cj^{-1}\ck\nabla_kH_0=\ck\nabla_k\tilde{H}_1,
$$
respectively. It then follows via an inductive argument, that $ \cR^n\cj\nabla_kH_0 = \ck \nabla_k
\tilde{H}_{n-1}=\cj\nabla_k\tilde{H}_n$, for $n>0$, where $\tilde{H}_0=H_0$.
\epf

\begin{cor}\label{COR:cflows}
(Validity of Item (i).) Let $\cj$ and $\ck$ be skew operators on the tangent bundle to $\cc$, and let $\{H_j\}$ be a sequence of functionals on $\cc$ indexed by $j$ in ${\mathbb Z}$.  If $\cj\nabla H_{k+1}=\ck\nabla H_k,$ then the $H_j$ pairwise commute, that is $\{H_j,H_k\}_\cj=0$.
\end{cor}
\pf  Simply using the skew-symmetry of $\ck$ and $\cj$, we have:
\begin{align}
\{H_j,H_k\}_\cj &= \left< \nabla H_j, \cj \nabla H_k \right>= \left< \nabla H_j, \ck \nabla H_{k-1} \right>= -\left< \ck \nabla H_j, \nabla H_{k-1} \right>\nonumber \\
               &= -\left< \cj \nabla H_{j+1}, \nabla H_{k-1} \right>= \left< \nabla H_{j+1}, \cj \nabla H_{k-1}
               \right>=\{H_{j+1},H_{k-1}\}_\cj. \label{CommFlow1}
\end{align}
Assume now $k>j$.  Then after $(k-j)$ iterations of the procedure yielding \eqref{CommFlow1}, we obtain $\{H_j,H_k\}_\cj = \{H_k,H_j\}_\cj,$ which implies $\{H_j,H_k\}_\cj=0$, as the $\cj$-bracket is skew.
\epf
\bigskip
\bigskip

\noindent{\bf Acknowledgments.} G. I. Lozano would like to thank Hermann Flaschka for useful
discussions and feedback.  N. M. Ercolani and G. I. Lozano were supported in part by NSF grant no. 0073087.
\bigskip

\section*{Appendix I: Key elements in the proof of Theorem \ref{THM:grads}}\label{APP:I}

Let $(\qdt_k, \rdt_k)^T=(a_k, b_k)^T$ denote an arbitrary variation of $(q_k, r_k)$, and $\nu_k$ denote a (vector) solution to the eigenvalue problem \eqref{ZS-disc}. The induced variation on
$\nu_k$ satisfies
\begin{equation}
\nut_{k+1}-\bpm z & q_k \\
r_k & 1/z\epm\nut_k=\bpm 0 & \qdt_k \\ \rdt_k & 0\epm\nu_k.
\label{nuVariation}
\end{equation}
Setting $\nut_k=\Phi(k)u_k$ for $\Phi(k)$ as in \eqref{PhiAsym}, and using the standard variation of
constants method to solve \eqref{nuVariation} for $\nut_k$ yields:
\begin{equation}
\Phi(k+1)u_{k+1}-\bpm z & q_k \\
r_k & 1/z\epm\Phi(k)u_k= \Phi(k+1)(\eu-1)u_k=\bpm 0 & \qdt_k \\ \rdt_k &0\epm\nu_k.
\end{equation}
Upon inverting $(\eu-1)$, one obtains
\begin{equation}
u_n=-\sum_{k=n}^\infty \Phi(k+1)^{-1}\bpm 0 & \qdt_k \\ \rdt_k &
0\epm\nu_k + c,\ \ \ \nut_n=\Phi(n)u_n. \label{ii}
\end{equation}
In particular, choosing $\nu_n=\psi_n$, we see that as $n\rightarrow\infty$, $\psi_n\sim(0,z^{-n})^T$.  So $\dot{\psi}_n$ vanishes and, clearly, so does the sum in \eqref{ii}, which implies $c=0$.  Multiplying the
resulting equation through by $z^n$ and focusing on the variation of just the second component of the resulting vector (i.e.,  $z^n \sit_n$) we see that:
\begin{align}
z^n \dot{\sit_n} &= -\sum_{k=n}^\infty \left[(z^n \Phi(n))\Phi(k+1)^{-1}\bpm 0 & \qdt_k \\
\rdt_k & 0\epm\psi_k\right]^{(2)}\label{psi2z^nVariation}\\
&= -\sum_{k=n}^\infty \left(z^n\fit_n\right)\left(\beta_{k+1}\rdt_k\sio_k+\alpha_{k+1}\qdt_k\sit_k\right)
+
\left(z^n\sit_n\right)\left(\delta_{k+1}\rdt_k\sio_k+\gamma_{k+1}\qdt_k\sit_k\right),
\notag
\end{align}
where
\begin{equation}
\bpm \alpha_{k+1} & \beta_{k+1} \\ \gamma_{k+1} & \delta_{k+1}\epm=\Phi(k+1)^{-1}=\frac{P(k+1)}{a}\bpm \sit_{k+1} & -\sio_{k+1} \\ -\fit_{k+1} &
\fio_{k+1}\epm.\label{PhiInv}
\end{equation}
(Here, one uses relations \eqref{phiScattRel} to determine $\det\Phi(k)=\fio_k\sit_k-\sio_k\fit_k=a/P(k)$ \cite{guada}\cite{APT}.)

Taking the limit of \eqref{psi2z^nVariation} as $n\rightarrow-\infty$ and using the asymptotics \eqref{PhiAsym} together with \eqref{PhiInv} above, we get:
\begin{equation}
\dot{\left(a/C_0\right)}=\sum_{k=-\infty}^\infty(P(k+1)/C_0)\left(-\fio_{k+1}\sio_k\ \rdt_k+\fit_{k+1}\sit_k\ \qdt_k\right).\label{aCoVariation!}
\end{equation}
Multiplying \eqref{aCoVariation!} by $C_0/a$, we obtain \eqref{aCoFinalVariation} and, hence, \eqref{yeah!}.

Focusing on the variation of $z^{-n}\sioh_n$ and arguing in a similar manner, one arrives at the identity
\begin{equation}
\stackrel{.}{\left({\hat{a}}/C_0\right)}=\sum_{k=-\infty}^\infty
(P(k+1)/C_0)\left(\fioh_{k+1}\sioh_k\
\rdt_k-\fith_{k+1}\sith_k\ \qdt_k\right),\label{ahCoVariation!}
\end{equation}
Theorem \ref{THM:grads} follows.

\section*{Appendix II: Key elements in the proofs of Theorem \ref{THM:resol} and Theorem \ref{THM:kc}}\label{APP:II}

We begin by defining the following matrices of squared eigenfunctions:
\begin{align}
V_k&=\Phi_k\bpm 0 & 0 \\ 0 & -1\epm \Phi_k^{-1}=\bpm \sio_k\fit_k
& -\sio_k\fio_k \\ \sit_k\fit_k & -\sit_k\fio_k \epm
\frac{P(k)}{a}=\bpm  C_k & -A_k\\ B_k & -D_k\epm, \label{VkMat} \\
W_k&=\Phi_k\bpm 0 & 0 \\0 & -1\epm
\Phi_{k+1}^{-1}=\bpm \sio_k\fit_{k+1}& -\sio_k\fio_{k+1} \\
\sit_k\fit_{k+1} & -\sit_k\fio_{k+1} \epm
\frac{P(k+1)}{a}=\bpm \ggk & -\aak\\ \bbk & -\ddk\epm.\label{WkMat}
\end{align}
Note that $V_k$'s entries are standard squared eigenfunctions, whereas $W_k$ is composed of semi-shifted products.  Observe also that the constant matrix defining these identities is chosen so
that the off-diagonal entries of $W_k$ are precisely the components of the vector of squared eigenfunctions defining the generating gradient $\delta_k\log a+\delta_k H_0$.

One can directly verify the following identities for $V_k$ and $W_k$:
\begin{align}
&\bullet \  W_{k+1}=\ce_k W_k \ce_{k+1}^{-1}; \label{1} \\
&\bullet \  W_k=V_k\Phi_k\Phi_{k+1}^{-1}=V_k \ce_k^{-1};
\label{2} \\
&\bullet \  (r_k\aak+z^{-1}\ddk)-(q_k\bbk+z\ggk)=1, \label{3}
\end{align}
where $\ce_k$ is, recall, the coefficient matrix for our eigenvalue problem \eqref{ZS-disc}, and the entries of $W_k$ are now written in terms of $\aak$, $\bbk$, $\ggk$ and $\ddk$, as in \eqref{WkMat}.

By writing \eqref{1} in terms of the new form of $W_k$ and tracking only the off-diagonal entries of the resulting matrix, one observes that
\begin{equation}
\eu\left[(1-r_kq_k) \abv \right]=\bpm z^{-2} & 0\\0 & z^2 \epm
\abv +\bpm r_kr_{k+1}\aak+z^{-1}(r_k\ggk+r_{k+1}\ddk) \\
          -q_kq_{k+1}\bbk-z(q_k\ddk+q_{k+1}\ggk)\epm. \label{4}
\end{equation}
Adding and subtracting the expression $\diag(z^{-2},z^2) [-r_kq_k(\bbk,-\aak)^T-(r_k,q_k)^T]$ on the right-hand-side of \eqref{4}, and using identity \eqref{3}, one may re-write it as
\begin{align}
\eu\left[(1-r_kq_k) \abv \right]&=\bpm z^{-2} & 0\\0 & z^2
\epm\left[(1-r_kq_k)\abv-\bpm
r_k\\q_k\epm\right] \label{+sh} \\
&+\bpm {r_{k+1}+z^{-2}r_k} &0 \\ 0& {q_{k+1}+z^2q_k}\epm\bpm q_k\bbk+z\ggk+1 \\
          -(r_k\aak+z^{-1}\ddk-1)\epm.\notag
\end{align}
By multiplying \eqref{+sh} by $\diag(z^2, z^{-2})\ed$, we obtain an almost perfectly symmetric equation in terms of the opposite shift:
\begin{align}
\ed\left[(1-r_kq_k) \abv \right]&=\bpm z^2 & 0\\0 & z^{-2}
\epm\left[(1-r_kq_k)\abv-\bpm
r_k\\q_k\epm\right]\label{-sh} \\
&+\bpm {r_{k-1}+z^2r_k} &0 \\ 0& {q_{k-1}+z^{-2}q_k}\epm\ed\bpm -(q_k\bbk+z\ggk) \\
          r_k\aak+z^{-1}\ddk\epm.\notag
\end{align}

Spelling out the relationship given in \eqref{2} one discovers:
\begin{equation}
D_k=z^{-1}\ddk-q_k\bbk, \ \ \ \mbox{ and } \ \ \ -C_k=r_k\aak-z\ggk, \label{VWrelations}
\end{equation}
so that \eqref{3} may be written as: $D_k-C_k=(z^{-1}\ddk-q_k\bbk)-(r_k\aak-z\ggk)=1$.  Using this
fact together with the previous relations \eqref{VWrelations}, we see that the $\delta\gamma$-dependent factors of \eqref{+sh} and \eqref{-sh} may be written as
\begin{equation}
(q_k\bbk+r_k\aak+C_k)\bpm 1\\-1\epm+\bpm 1\\0\epm,  \ \ \ \mbox{ and } \ \ \
(q_k\bbk+r_k\aak+C_k) \bpm 1\\1\epm+\bpm 0\\1\epm,\label{cpiece}
\end{equation}
respectively.  Finally, writing out $V_{k+1}$ in terms of its entries, we find that $C_k$ satisfies a difference equation given in terms of $\aak$ and $\bbk$, namely
\begin{equation}
(\eu -1)C_k=q_k\bbk+r_k\aak.\nonumber
\end{equation}
Taking $C_k=-\sum_{j=k}^\infty q_j\beta_j+r_j\alpha_j$, and substituting this expression into Eqs. \eqref{cpiece}, we obtain
\begin{align}
\bullet &\ \ \  \bpm {q_k\beta_k+z\gamma_k+1} \\ {-(r_k\alpha_k+z^{-1}\delta_k-1)}\epm  =
\bpm  -J^+_{k+1}(q_j) & J^+_{k+1}(r_j) \\ J^+_{k+1}(q_j) & -J^+_{k+1}(r_j) \epm\abv +\bpm 1 \\ 0\epm, \notag\\
\bullet & \ \ \  \bpm {q_k\beta_k+z\gamma_k} \\ {r_k\alpha_k+z^{-1}\delta_k}\epm  =
\bpm  -J^+_{k+1}(q_j) & J^+_{k+1}(r_j) \\ -J^+_{k+1}(q_j) & J^+_{k+1}(r_j) \epm\abv +\bpm 0 \\ 1\epm, \label{abgd2}
\end{align}
expressing the $\delta\gamma$-dependent factors in \eqref{+sh} and \eqref{-sh} in terms of the sum operator $\jkp$ in the formula for $L$, and the non-standard squared eigenfunctions $\aak$, $\bbk$.

Equations \eqref{+sh}, \eqref{-sh} used in conjunction with \eqref{abgd2} (as well as \eqref{+shHat}, \eqref{-shHat} together with \eqref{abgd2Hat}) play key roles in the proofs of Theorem \ref{THM:resol} and Theorem \ref{THM:kc}, as indicated therein.  The latter triple arises by considering the matrix of (semi-shifted) squared eigenfunctions tied to the generating function $\abvhr$ described in \eqref{ahGrads-psi1}, namely,
\begin{equation}
\bpm -\gammah_k & \alphah_k \\ -\betah_k & \deltah_k \epm = \frac{P(k+1)}{\ah}
\bpm -\sioh_k\fith_{k+1} & \sioh_k\fioh_{k+1} \\ -\sith_k\fith_{k+1} & \sith_k\fioh_{k+1}\epm, \ \ (r_k\alphah_k+z^{-1}\deltah_k)-(q_k\betah_k+z\gammah_k)=1, \label{WkMatHat}
\end{equation}
Starting with \eqref{WkMatHat}, and applying a procedure analogous to the one described in this appendix, we obtain:
\begin{align}
\eu\left[(1-r_kq_k) \abvh \right]=&\bpm z^{-2} & 0\\0 & z^2
\epm\left[(1-r_kq_k)\abvh+\bpm
r_k\\q_k\epm\right]\notag\\
-&\bpm {r_{k+1}+z^{-2}r_k} & 0 \\ 0& {q_{k+1}+z^2q_k}\epm\bpm q_k\betah_k+z\gammah_k+1 \\
          -(r_k\alphah_k+z^{-1}\deltah_k-1)\epm, \label{+shHat}
\end{align}
\begin{align}
\ed\left[(1-r_kq_k) \abvh \right]&=\bpm z^2 & 0\\0 & z^{-2}
\epm\left[(1-r_kq_k)\abvh+\bpm
r_k\\q_k\epm\right]\notag \\
&+\bpm {r_{k-1}+z^2r_k} &0 \\ 0& {q_{k-1}+z^{-2}q_k}\epm\ed\bpm q_k\betah_k+z\gammah_k \\
          -(r_k\alphah_k+z^{-1}\deltah_k)\epm,\label{-shHat}
\end{align}
where
\begin{equation}
\bpm {q_k\betah_k+z\gammah_k} \\
{-(r_k\alphah_k+z^{-1}\deltah_k)+1}\epm=\bpm  J^+_{k+1}(q_j) & -J^+_{k+1}(r_j)
\\ -J^+_{k+1}(q_j) & J^+_{k+1}(r_j) \epm\abvh. \label{abgd2Hat}
\end{equation}
as just mentioned above.

\section{Conclusion}

To summarize, this paper shows  that the AL hierarchy  can be explicitly viewed as a hierarchy of commuting flows which: (a) are Hamiltonian with respect to both the standard, local Poisson operator $\cj$, and a new non-local, skew, almost Poisson operator $\ck$, on the appropriate space; (b) can be recursively generated from the {\it recursion} operator $\cR=\ck\cj^{-1}$.  In addition, the proof of 
these facts relies upon two new pivotal resolvent identities which suggest a general method for
uncovering bi-Hamiltonian structures for other families of discrete, integrable equations.

Another result stemming from the current research is the clarification of the geometric framework that underlies a certain class of geodesic linkages evolving on the sphere \cite{DS, guada}.   A linkage on a Riemannian manifold is essentially defined by specifying a sequence of points connected by geodesic arcs.  A closed linkage is usually called a polygon.  Such linkages  are related to the AL hierarchy via the evolution for their ``discrete" geodesic curvature \cite{DS}.  In this regard, Lozano's preliminary results include a geometric interpretation of a compatibility condition associated to a Lax pair for AL and, consequently, a bijective correspondence between discrete, integrable mKdV flows (also AL flows) and linkage flows. (For details on this, see \cite{guada}; also see \cite{DS, Langer, LanPerKdV} for background in terms of continuous-analogs of the linkage models.)

Let us now close by summarizing some of the many possible avenues for further research.  First off,
a definite answer to the question of whether or not $\ck$ actually defines a Poisson structure would be
desirable  on several counts.  If $\ck$ were Poisson, one could turn to exploring the possible connections between the $\ck$-induced Poisson bracket and the bi-Hamiltonian
structure for finite AL described by Faybusovich and Gekhtman \cite{FGekht}.  This work considers the AL hierarchy within the larger class of full Toda flows in $sl(n)$
and presents a bi-Hamiltonian formulation for (finite) AL stemming from the bi-Hamiltonian structure
of these Toda flows.  

On the other hand, an obstruction to the Jacobi identity would place the $\ck$-bracket in the category
of almost Poisson structures and could perhaps steer the investigation in the direction of non-holonomic mechanical systems \cite{Bloch}.  Such systems possess an underlying Hamiltonian structure that is 
(strictly) almost Poisson. 

The elucidation of new bi-Hamiltonian structures and their connection with the evolution of non-stretching classes of linkages could also be pursued in the context of the de-focusing Ablowitz--Ladik system (obtained from \eqref{ALrq} in the reduction $r_k=\overline{q}_k$) and other discrete integrable equations.  This would serve as a test of the robutness of our methods for deriving operators such
as $\ck$, $\cR$, and resolvent identities such as those obtained for $L$ and $L^{-1}$.  It would be 
good, for instance, to understand how recursion operators (such as $\cR$,
or even $L$ and $L^{-1}$) are encoded in the squared eigenfunctions of a linear problem (through
resolvent identities of the right kind.)

The connection with ``physical linkage" spaces could also be pursued further, both for discrete,
integrable mKdV, AL and potentially for other discrete integrable equations.  In the context of discrete
integrable mKdV, AL and non-stretching spherical linkages for instance, one should aim at
understanding the linkage recursion schemes proposed in \cite{guada} in Poisson-geometric terms.
Ideally, a well-defined lift of $\cR$ to the space of non-stretching linkages could be defined and then parsed out as the composition of two Poisson (or perhaps one Poisson and one almost Poisson) operators.    Connections with results of Langer and Perline (for the case of continuous NLS and the FM model \cite{LanPer},) and  Kapovich and Millson (regarding the symplectic geometry of non-stretching polygons\cite{KM, KMhyp}, could also be explored and addressed in the appropriate context.

\end{document}